\documentclass[11pt]{article}
\usepackage{amsmath}
\usepackage{amsfonts}
\usepackage{amssymb}
\usepackage{epsfig}
\usepackage[francais,english]{babel}

\newtheorem{thm}{Theorem}[section]
\newtheorem{lem}[thm]{Lemma}
\newtheorem{corollary}[thm]{Corollary}
\newtheorem{proposition}[thm]{Proposition}
\newtheorem{conjecture}{Conjecture}
\newtheorem{remark}[thm]{Remark}

\newcommand{\A}{{\cal A}}
\newcommand{\B}{{\cal B}}
\newcommand{\C}{{\mathbb C}}

\newcommand{\fs}{{\sf f}}
\newcommand{\Fs}{{\sf F}}
\newcommand{\G}{{\cal G}}
\def\II{{\mathbb I}}
\newcommand{\Ks}{{\sf K}}
\newcommand{\ks}{{\sf k}}
\newcommand{\LL}{{\cal L}}
\newcommand{\Ls}{{\sf L}}
\newcommand{\ls}{{\sf l}}
\newcommand{\M}{{\cal M}}
\newcommand{\NN}{{\mathbb N}}
\newcommand{\Ns}{{\sf N}}
\newcommand{\Ps}{{\sf P}}
\newcommand{\Q}{{\cal Q}}
\newcommand{\Qs}{{\sf Q}}

\newcommand{\Rs}{{\sf R}}
\newcommand{\RR}{{\mathbb R}}
\newcommand{\Ss}{{\sf S}}
\def\Perm{{\cal S}}
\newcommand{\TT}{{\mathbb T}}
\newcommand{\ZZ}{{\mathbb Z}}

\def\jump{\vspace{-0.2cm}\\}

\def\G1bar{\bar{G}_1}

\def\L1bar{\bar{L}_1}

\def\S1bar{\bar{S}_1}
\def\Snbar{\bar{S}_n}

\newenvironment{proof}{
    \noindent{\it Proof:} \hspace*{1em}}{
    \hspace*{\fill} $\square$\bigskip \jump}

\begin{document}
\title{Quadri-tilings of the plane}
\author{B\'eatrice de Tili\`ere\vspace{0.5cm}
\thanks{{\small\texttt{beatrice.detiliere@math.unizh.ch}}}\\
{\small Institut f\"ur Mathematik,
       Universit\"at Z\"urich,
       Winterthurerstrasse 190,
       CH-8057 Z\"urich.}
}
\date{}
\maketitle

\begin{abstract}
We introduce {\em quadri-tilings} and show that they are in
bijection with dimer models on a {\em family} of graphs $\{R^*\}$
arising from rhombus tilings. Using two height functions, we
interpret a sub-family of all quadri-tilings, called {\em
triangular quadri-tilings}, as an interface model in dimension
$2+2$. Assigning ``critical" weights to edges of $R^*$, we prove
an explicit expression, only depending on the local geometry of
the graph $R^*$, for the minimal free energy per fundamental
domain Gibbs measure; this solves a conjecture of \cite{Kenyon1}.
We also show that when edges of $R^*$ are asymptotically far
apart, the probability of their occurrence only depends on this
set of edges. Finally, we give an expression for a Gibbs measure
on the set of {\em all} triangular quadri-tilings whose marginals
are the above Gibbs measures, and conjecture it to be that of
minimal free energy per fundamental domain.
\end{abstract}

\section{Introduction}

In this paper, we introduce {\bf quadri-tilings} and a sub-family
of all quadri-tilings called {\bf triangular quadri-tilings}. In order
to explain the originality of this model, let us go back to the yet
classical domino and lozenge tilings, see for example
\cite{EKLP,Kast1,Kenyon,Thurston}. 
Both models are dimer models (see below) on a
{\em fixed} graph, the square lattice $\ZZ^2$, and the equilateral
triangular lattice $\TT$, respectively. By the means of a {\em height
function}, they can be interpreted as random discrete interfaces in dimension
$2+1$, that is as random discrete surfaces of dimension $2$ in a space of
dimension $3$ that have been projected to the plane
\cite{Thurston}. Keeping this in mind, one can now explain the
interesting feature of quadri-tilings: they correspond to
dimer models on a {\em family} of graphs, instead of a fixed graph as
was the case up to now. Moreover, by the means of
{\em two} height functions, triangular quadri-tilings can be
interpreted as random interfaces in dimension $2+2$, that is random
surfaces of dimension $2$ in a space of dimension $4$, that have been
projected to the plane. It is the first such model arising from dimer
models. Using tools of the dimer model, we are able to give an
explicit expression for a Gibbs measure on triangular quadri-tilings,
as well as a surprising property of the asymptotics of this
measure. In the course of doing so, we prove a conjecture of \cite{Kenyon1}.
\jump

\begin{figure}[h]
\begin{center}
\includegraphics[width=\linewidth]{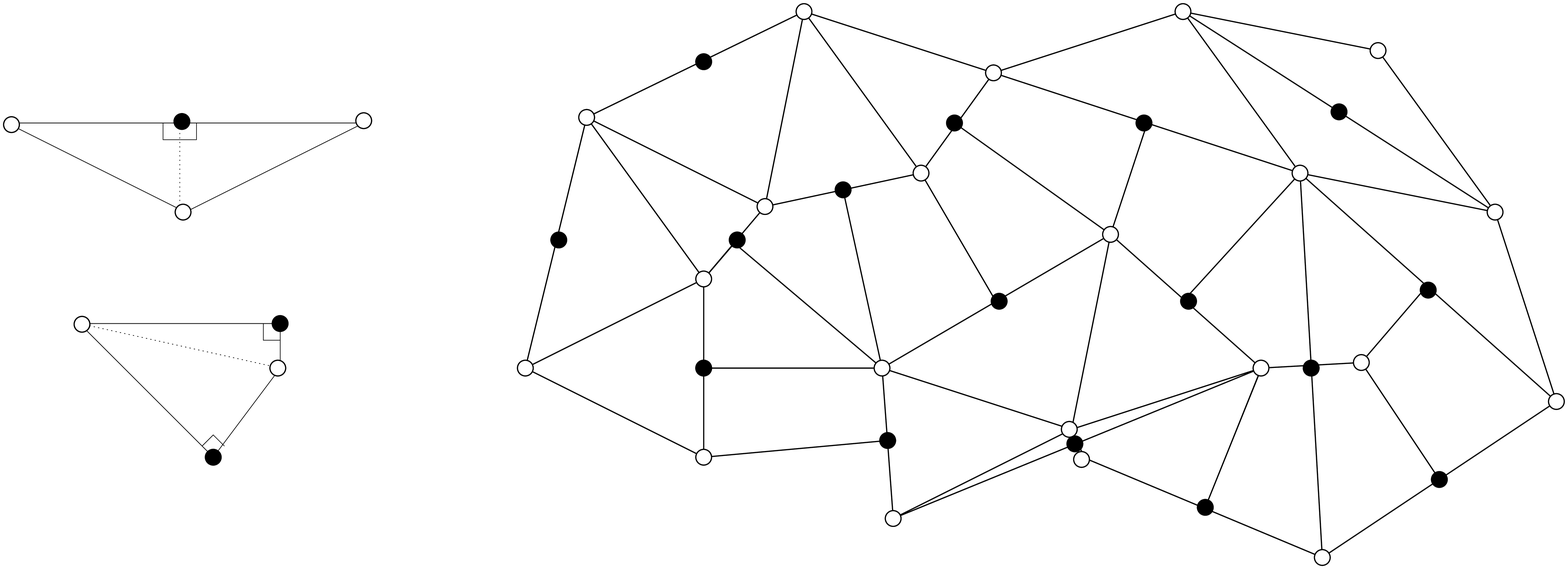}
\end{center}
\caption{Two types of quadri-tiles (left), and a quadri-tiling
(right).}\label{fig1}
\end{figure}

\noindent Quadri-tilings are defined as follows. Consider the set
of right triangles whose hypotenuses have length two. Color the
vertex at the right angle black, and the other two vertices white.
A {\bf quadri-tile} is a quadrilateral obtained from two such
triangles in two different ways: either glue them along the
hypotenuse, or supposing they have a leg of the same length, glue
them along this edge matching the black (white) vertex to the
black (white) one, see Figure \ref{fig1}. Note that both types of
quadri-tiles have four vertices. A {\bf quadri-tiling} of the
plane is an edge-to-edge tiling of the plane by quadri-tiles that
respects the coloring of the vertices, that is black (resp. white)
vertices are matched to black (resp. white) ones. An example of
quadri-tiling is given in Figure \ref{fig1}. In all that follows,
we consider quadri-tilings of the plane that use finitely many
different quadri-tiles up to isometry.

\noindent The goal of Section \ref{sec2} is to precisely describe the
features of quadri-tilings. In order to give some insight, let us
define {\bf $2$-tiling models} or equivalently {\bf dimer models}.
A {\bf $2$-tile} of an infinite graph $G$ is a polygon made of two
adjacent inner faces of $G$, and a {\bf $2$-tiling} of $G$ is a
covering of $G$ with $2$-tiles, such that there are no holes and
no overlaps. The {\bf dual graph} $G^*$ of $G$ is the graph whose
vertices correspond to faces of $G$, two vertices of $G^*$ being
joined by an edge if the corresponding faces are adjacent. A {\bf
dimer configuration} of $G^*$, also called {\bf perfect matching},
is a subset of edges of $G^*$ which covers each vertex exactly
once. Then $2$-tilings of the graph $G$ are in bijection with
dimer configurations of the dual graph $G^*$, as explained by the
following correspondence. Denote by $f^*$ the dual vertex of a
face $G$, and consider an edge $f^* g^*$ of $G^*$. We say that the
$2$-tile of $G$ made of the adjacent faces $f$ and $g$ is the
$2$-tile {\bf corresponding} to the edge $f^* g^*$. Then,
$2$-tiles corresponding to edges of a dimer configuration of $G^*$
form a $2$-tiling of $G$. Let us denote by $\M(G^*)$ the set of
perfect matching of the graph $G^*$.\jump 

\noindent Prior to describing Section
\ref{sec2}, we need one more definition. If $\Rs$ is a rhombus
tiling of the plane, then the corresponding {\bf
rhombus-with-diagonals tiling}, denoted by $R$, is the graph
obtained from $\Rs$ by adding the diagonals of the rhombi. In the whole of this paper, we
suppose that rhombus tilings of the plane have only finitely many
rhombus angles.\jump

\noindent In Section \ref{subsec21}, we prove the main feature of
quadri-tilings, i.e. that they correspond to $2$-tilings on a family of
graphs, which consists of rhombus-with-diagonals tilings. More
precisely, if $T$ is a quadri-tiling, then by a geometric construction,
we associate to $T$ a rhombus-with-diagonals tiling $R(T)$, such that
$T$ is a $2$-tiling of $R(T)$. The corresponding rhombus tiling
$\Rs(T)$ is called the {\bf underlying rhombus tiling} of $T$.\jump

\begin{figure}[h]
\begin{center}
\includegraphics[width=\linewidth]{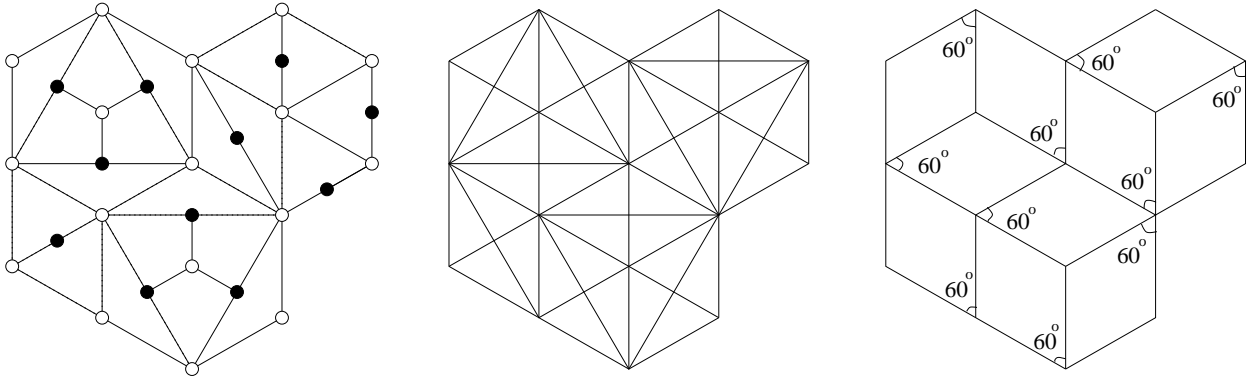}
\end{center}
\caption{Triangular quadri-tiling $T$ (left), underlying
lozenge-with-diagonals tiling $L(T)$ (middle), corresponding lozenge
tiling $\Ls(T)$ (right).}\label{fig15}
\end{figure}

\noindent {\bf Triangular quadri-tilings} consist in the
sub-family of all quadri-tilings whose underlying tiling is a lozenge tiling, where 
{\bf lozenges} are defined to be $60^{\circ}$-rhombi; refer to Figure
\ref{fig15} for an example, and to Section \ref{subsec21} for the construction of
the underlying tiling. In
order to distinguish general rhombus tilings, denoted $\Rs$, from lozenge tilings,
we denote the latter by $\Ls$. The set of all
triangular quadri-tilings up to isometry is denoted by $\Q$. Note that
$\Q$ corresponds to two superposed dimer models. Indeed, let $T$ be a
triangular quadri-tiling, then $T$ is a
$2$-tiling of its underlying lozenge-with-diagonals tiling $L(T)$,
moreover the corresponding lozenge tiling $\Ls(T)$ is a $2$-tiling
of the equilateral triangular lattice $\TT$.

\noindent Section \ref{subsec22} consists in the geometric
interpretation of triangular quadri-tilings using height functions. On the
vertices of every triangular quadri-tiling $T$, we define a $\ZZ$-valued
function $h_1$, called the {\bf first height function},
corresponding to the ``height'' of $T$ interpreted as a $2$-tiling
of its underlying lozenge-with-diagonals tiling $L(T)$. Then, we
assign a {\bf second height function} $h_2$
(Thurston's height function on lozenges \cite{Thurston})
corresponding to the height of $\Ls(T)$ interpreted as a
$2$-tiling of $\TT$, see Figure \ref{fig3}. Hence triangular quadri-tilings are
characterized by two height functions, and so can be interpreted
as discrete interfaces in dimension $2+2$. In Section \ref{subsec23}, we
give elementary operations that allow to transform any
triangular quadri-tiling of a simply connected region into any other.\jump

\noindent The dimer model belongs to the field of statistical
mechanics, hence there are natural measures to consider, Boltzmann
and Gibbs measures, which are defined as follows. Let $G$ be an
infinite graph, and let $\nu$ be a positive weight function on the
edges of $G^*$. Consider a finite sub-graph $G^1$ of $G$, then the
{\bf Boltzmann measure} on the set of dimer configurations
$\M({G^1}^*)$ of ${G^1}^*$, corresponding to the weight function
$\nu$, is defined by
\begin{equation*}\label{boltzman1}
\mu^1(M)=\frac{ \prod_{e \in M} \nu(e)}{Z({G^1}^*,\nu)},
\end{equation*}
where $Z({G^1}^*,\nu)= \sum_{M \in \M({G^1}^*)} \prod_{e \in M}
\nu(e)$ is the {\bf dimer partition function}. A {\bf Gibbs
measure} is a probability measure on $\M(G^*)$ with the following
property. If the matching in an annular region is fixed, the
matchings inside and outside of the annulus are independent of
each other, and the probability of any interior matching $M$ is
proportional to $\prod_{e \in M} \nu(e)$.\jump

\noindent An important question in solving a dimer model is the
study of local statistics, i.e. to obtain an explicit expression
for the set of Gibbs measures. Kenyon, Okounkov and Sheffield
\cite{KeOS} give such an expression for the two-parameter family
of Gibbs measures on dimer configurations of doubly periodic
bipartite graphs. The expression they obtain involves the limiting
inverse Kasteleyn matrix which is hard to evaluate in general,
often implying elliptic integrals. In another paper
\cite{Kenyon1}, for graphs $G$ which have bipartite duals and
satisfy a geometric condition called {\bf isoradiality}, Kenyon
defines a specific weight function on the edges of $G^*$ called
the {\bf critical} weight function. He also defines the Dirac
operator $K$ indexed by the vertices of $G^*$, and gives an
explicit expression for its inverse $K^{-1}$ (see also Sections
\ref{subsec31} and \ref{subsec32}). The expression for $K^{-1}$
has the interesting property of only depending on the {\em local
geometry} of the graph. Kenyon conjectures that $K^{-1}$ is, in
some sense to be determined, the limiting inverse Kasteleyn
matrix. In Section \ref{sec3}, we consider a general
rhombus-with-diagonals tiling of the plane $R$. It has the
property of being an isoradial graph, so that we assign the
critical weight function to edges of $R^*$. Theorem \ref{thm2} of
Section \ref{subsec33} (see also Theorem \ref{thm0} below) proves
an explicit expression for a Gibbs measure $\mu^R$ on $\M(R^*)$,
as a function of $K$ and $K^{-1}$.\jump

\noindent For every subset of edges $e_1=w_1 b_1,\ldots,e_k=w_k
b_k$ of $R^*$, the {\bf cylinder} $\{e_1,\ldots,e_k\}$ is defined
to be the set of dimer configurations of $R^*$ that contain these
edges. Then
\begin{thm}\label{thm0}
There is a probability measure $\mu^R$ on $\M({R}^*)$ such that,
for every cylinders $\{e_1,\ldots,e_k\}$ of $R^*$,
\begin{equation}\label{3}
\mu^R(e_1,\ldots,e_k)=\left( \prod_{i=1}^{k} K(w_i,b_i)\right)
\det_{1 \leq i,\;j \leq k} \left(K^{-1}(b_i,w_j)\right).
\end{equation}
Moreover $\mu^R$ is a Gibbs measure on $\M(R^*)$. When $R^*$ is
doubly periodic, $\mu^R$ is the unique Gibbs measure which has
minimal free energy per fundamental domain among the two-parameter
family of ergodic Gibbs measures of \cite{KeOS}.
\end{thm}
\begin{enumerate}
    \item[$\bullet$] Note that we do not ask the graph $R^*$ to be
periodic. The proof of Theorem \ref{thm0} is the subject of
Section \ref{sec4}. The argument in the case where $R^*$ is not
periodic relies on the argument in the case where $R^*$ is doubly
periodic, combined with a non-trivial geometric property of
rhombus tilings proved in Proposition~\ref{prop1}: ``every finite
simply connected sub-graph of a rhombus tiling can be embedded in
a periodic rhombus tiling of the plane."
    \item[$\bullet$] From the proof of Theorem
\ref{thm0}, it appears that the statement is true for all doubly
periodic isoradial graphs with bipartite duals. Hence,
Theorem~\ref{thm0} solves the conjecture of \cite{Kenyon1} of
interpreting the inverse Dirac operator as the limiting inverse
Kasteleyn matrix. Moreover, the fact that the measure $\mu^R$ is
of minimal free energy per fundamental domain makes it of special
interest among the two-parameter family of ergodic Gibbs measures
of \cite{KeOS}.
    \item[$\bullet$] Using the locality property of $K^{-1}$ mentioned
above, we deduce that the expression (\ref{3}) only depends on the
local geometry of the graph, hence it yields an easy way of
computing local statistics explicitly. This is very surprising in
regards of the expression obtained in \cite{KeOS}, and we believe
this locality property to be true only in the isoradial case with
critical weights.
\end{enumerate}
In Section \ref{sec5}, we extend the notion of Gibbs measure to
the set of all triangular quadri-tilings $\Q$. Then, as a
corollary to Theorem \ref{thm0} we deduce an explicit expression
for such a Gibbs measure $\mu$, and conjecture it to be that of
minimal free energy per fundamental domain among a four parameter
family of Gibbs measures.\jump

\noindent In Section \ref{sec6}, we consider a general
rhombus-with-diagonals tiling of the plane $R$. We assign the
critical weight function to edges of $R^*$, and let $K$ be the
Dirac operator indexed by vertices of $R^*$. Theorem \ref{thm7} of
Section \ref{subsec61} (see also Theorem \ref{thm1} below)
establishes that asymptotically (as $|b-w|\rightarrow\infty$) and
up to the second order term, $K^{-1}(b,w)$ only depends on the
rhombi to which the vertices $b$ and $w$ belong, and else is
independent of the structure of the graph $R$. For a general
isoradial graph, Kenyon \cite{Kenyon1} gives an asymptotic formula
for $K^{-1}(b,w)$ which depends on the angles of an edge-path from
$w$ to $b$. Hence, it is an interesting and surprising fact that
the dependence on the edges along the path should asymptotically
disappear in the case of rhombus-with-diagonals tilings.
\begin{thm}\label{thm1}
As $|b-w|\rightarrow\infty$, $K^{-1}(b,w)$ is equal to\jump

\noindent {\small $ \frac{1}{2 \pi}\left(
\frac{1}{b-w}+\frac{e^{-i(\theta_1+\theta_2)}}{\bar{b}-\bar{w}}
\right) +\frac{1}{2 \pi}\left(\frac{e^{2
i\theta_1}+e^{2i\theta_2}}{(b-w)^3}+ \frac{e^{-i(3
\theta_1+\theta_2)}+e^{-i(\theta_1+3\theta_2)}}{(\bar{b}-\bar{w})^3}\right)
+O\left( \frac{1}{|b-w|^3}\right).$}
\end{thm}

\noindent As a consequence of Theorem \ref{thm1}, we deduce that
when edges $e_1,\ldots,e_k$ of $R^*$ are asymptotically far apart,
$\mu^R(e_1,\ldots,e_k)$ only depends on the rhombi to which the
edges $e_1,\ldots,e_k$ belong, and else is independent of the
structure of the graph $R$ (Corollary \ref{cor1}). We conclude by
giving a consequence of Corollary \ref{cor1} for the measure $\mu$
on triangular quadri-tilings (Corollary \ref{cor3}).\\

\noindent {\it Acknowledgments:} We thank Richard Kenyon for
proposing the quadri-tiling model and asking the questions related
to it; we are grateful to him for the many enlightening
discussions. We also thank the referee for the many pertinent
remarks which have helped to increase the quality of this paper.

\section{Features of quadri-tilings}\label{sec2}

\subsection{Underlying rhombus-with-diagonals tilings}\label{subsec21}

\begin{lem}\label{lem1}
Quadri-tilings are in one-to-one correspondence with $2$-tilings
of graphs which are rhombus-with-diagonals tilings of the plane.

\end{lem}
\begin{proof}
Consider a quadri-tiling of the plane $T$. Denote by $R$ the
tiling of the plane obtained from $T$ by drawing, for each
quadri-tile, the edge separating the two right triangles. Let $b$
be a black vertex of $R$, denote by $w_1,\ldots,w_k$ the neighbors
of $b$ in cclw (counterclockwise) order. In each right triangle,
the black vertex is adjacent to two white vertices, and since the
gluing respects the coloring of the vertices, $w_1,\ldots,w_k$ are
white vertices. Moreover, $b$ is at the right angle, so $k=4$ and
the edges $w_1 w_2$, $w_2 w_3$, $w_3 w_4$, $w_4w_1$ are
hypotenuses of right triangles. Therefore $w_1,\ldots,w_4$ form a
side-length-$2$ rhombus, and $b$ stands at the crossing of its
diagonals. This is true for any black vertex $b$ of $R$, so $R$ is
a rhombus-with-diagonals tiling of the plane, and $T$ is a
$2$-tiling of $R$.
\end{proof}
As a consequence of Lemma \ref{lem1}, a quadri-tiling $T$ is a
$2$-tiling of a unique rhombus-with-diagonals tiling, which we
call the {\bf underlying rhombus-with-diagonals tiling}, and
denote by $R(T)$, see Figure \ref{fig15}.

\subsection{Height functions}\label{subsec22}

\noindent We define a {\bf first height function} $h_1$ on
vertices of every quadri-tiling $T$. Moreover, when $T$ is a
triangular quadri-tiling, we define a {\bf second height function}
$h_2$ on vertices of $T$. Using $h_1$ and $h_2$, we interpret
triangular quadri-tilings as discrete $2$-dimensional surfaces in
a $4$-dimensional space projected to the plane.

\subsubsection{First height function}\label{subsubsec221}

\noindent Consider a quadri-tiling of the plane $T$, then $T$ is a
$2$-tiling of its underlying rhombus-with-diagonals tiling $R(T)$.
In order to define the first height function $h_1$, we need a
bipartite coloring of the faces of $R(T)$, which is given by the
following.

\begin{lem}\label{lem2}
Let $\Rs$ be a rhombus tiling of the plane, and $R$ be the
corresponding rhombus-with-diagonals tiling. Then $R$ has a
bipartite coloring of its faces which is also a bipartite coloring
of the vertices of $R^*$.
\end{lem}

\begin{proof}
Cycles corresponding to the faces of the graph $\Rs$ have length
four, thus $\Rs$ has a bipartite coloring of its vertices, say
black and white. Consider a face of $R$ and orient its boundary
edges cclw. If the white vertex of the hypotenuse-edge comes
before the black one, assign color black to the face, else assign
color white. This defines a bipartite coloring of the faces of
$R$, which is also a bipartite coloring of the vertices of $R^*$
(see Figure \ref{fig2}).\\
\begin{figure}[h]
\begin{center}
\includegraphics[height=3cm]{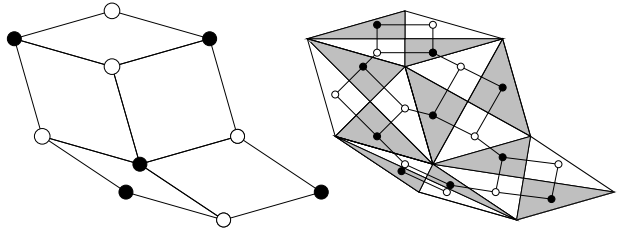}
\end{center}
\caption{Bipartite coloring of the vertices of $\Rs$ (left), and
corresponding bipartite coloring of the faces of $R$ and of the
vertices of ${{R}^*}$ (right).} \label{fig2}
\end{figure}
\end{proof}

\noindent Consider the bipartite coloring of the faces of $R(T)$. Orient the
edges around the black faces cclw, edges around the white faces
are then oriented cw, and define $h_1$ on the vertices of $T$ as
follows. Fix a vertex $v_1$ on a boundary edge of a rhombus of
$R(T)$, and set $h_1(v_1)=0$. For every other vertex $v$ of $T$,
take an edge-path $\gamma_1$ from $v_1$ to $v$ which follows the
boundaries of the quadri-tiles of $T$. The first height function
$h_1$ changes by $\pm 1$ along each edge of $\gamma_1$: if an edge
is oriented in the direction of the path, then $h_1$ increases by
$1$, if it is oriented in the opposite direction, then $h_1$
decreases by $1$. The value $h_1(v)$ is independent of the path
$\gamma_1$ because the plane is simply connected, and the height
change around any quadri-tile is zero. An example of computation
of $h_1$ is given in Figure \ref{fig3}.\jump

\noindent The following lemma gives a bijection between
$2$-tilings of a rhombus-with-diagonals tiling $R$ and first
height functions defined on vertices of $R$.

\begin{lem}\label{lem3}
Fix a vertex $v_1$ on a boundary edge of a rhombus of $R$. Let
$\tilde{h_1}$ be a $\ZZ$-valued function on the vertices of $R$
satisfying the following two conditions:
\begin{itemize}
\item[$\bullet$] $\tilde{h_1}(v_1)=0$, \item[$\bullet$]
$\tilde{h_1}(v)=\tilde{h_1}(u)+1$, or
$\tilde{h_1}(v)=\tilde{h_1}(u)-2$, for any edge $uv$ oriented from
$u$ to $v$.
\end{itemize}
Then, there is a bijection between functions $\tilde{h_1}$
satisfying these two conditions and $2$-tilings of $R$.
\end{lem}

\begin{proof}
The idea of the proof closely follows \cite{EKLP}. If $T$ is a
$2$-tiling of $R$, then the first height function defined above
satisfies the two conditions of the lemma: if an edge $u v$,
oriented from $u$ to $v$, belongs to the boundary of a
quadri-tile, it satisfies $h_1(v)=h_1(u)+1$, else if
it lies across a quadri-tile, it satisfies $h_1(v)=h_1(u)-2$.\\
Conversely, consider a $\ZZ$-valued function $\tilde{h_1}$ as in
the lemma. Then, anytime there is an edge $u v$ satisfying $|
\tilde{h_1}(v)-\tilde{h_1}(u) |=2$, put a quadri-tile made of the
two right triangles adjacent to this edge. This defines a
$2$-tiling of $R$.
\end{proof}

\subsubsection{Second height function}\label{subsec222}\ \jump

\noindent Consider a triangular quadri-tiling $T$. Let $L(T)$ be
its underlying lozenge-with-diagonals tiling, and $h_1$ be the
first height function on vertices of $T$. The lozenge tiling
$\Ls(T)$ corresponding to $L(T)$ is a $2$-tiling of the
equilateral triangular lattice $\TT$. Moreover $\TT$ has a
bipartite coloring of its faces, say black and white. Orient the
edges around the black faces cclw, edges around the white faces
are then oriented cw. Thurston \cite{Thurston} defines the second
height function $h_2$ as follows: choose a vertex $v_2$ of
$\Ls(T)$, and set $h_2(v_2)=0$. For every other vertex $v$ of
$\Ls(T)$, take an edge-path $\gamma_2$ from $v_2$ to $v$ which
follows the boundaries of the lozenges of $\Ls(T)$. The second
height function $h_2$ changes by $\pm 1$ along each edge of
$\gamma_2$: if an edge is oriented in the direction of the path,
then $h_2$ increases by $1$, if it is oriented in the opposite
direction, then $h_2$ decreases by $1$. The value $h_2(v)$ is
independent of the path $\gamma_2$. For
convenience, we choose $v_2$ to be the same vertex as $v_1$, and
denote this common vertex by $v_0$, so that $h_1(v_0)=h_2(v_0)=0$.
An analog to Lemma \ref{lem3} gives a bijection between second
height functions and lozenge tilings of the plane, hence we deduce
that triangular quadri-tilings are characterized by $h_1$ and
$h_2$.\jump

\begin{figure}[h]
\begin{center}
\includegraphics[height=5.6cm]{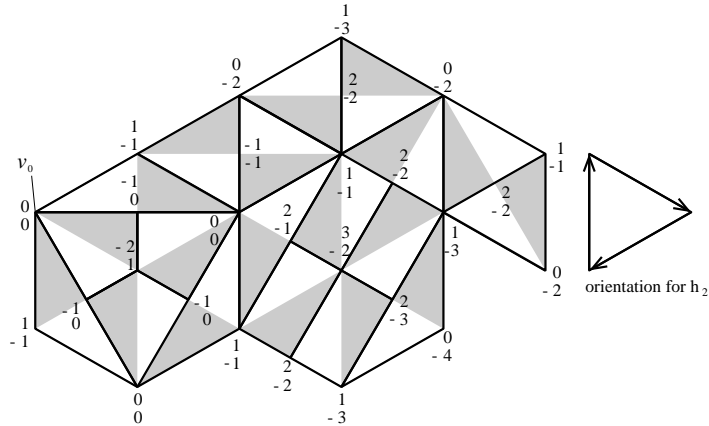}
\end{center}
\caption{Triangular quadri-tiling with height
functions $h_1$ (above) and $h_2$ (below).} \label{fig3}
\end{figure}

\noindent Let us define a natural value for the second height
function at the vertex in the center of the lozenges of $\Ls(T)$.
When going cclw around the vertices of a lozenge $\ell$ of
$\Ls(T)$, starting from the smallest value of $h_2$ say $h$,
vertices take on successive values $h, h+1, h+2, h+1$, so that we
assign value $h+1$ to the vertex in the center of the lozenge
$\ell$. An example of computation of $h_2$ is given in Figure
\ref{fig3}.

\noindent In Thurston's geometric interpretation \cite{Thurston},
a lozenge tiling is seen as a surface $S$ in $\ZZ^3$ (where the
diagonals of the cubes are orthogonal to the plane) that has been
projected orthogonally to the plane. The surface $S$ is determined
by the height function $h_2$. In a similar way, a triangular
quadri-tiling of the plane $T$ can be seen as a surface $S_1$ in a
4-dimensional space that has been projected to the plane; $S_1$
can also be projected to $\widetilde{\ZZ}^3$ ($\widetilde{\ZZ}^3$
is the space $\ZZ^3$ where cubes are drawn with diagonals on their
faces), and one obtains a surface $S_2$. When projected to the
plane, $S_2$ is the underlying lozenge-with-diagonals tiling
$L(T)$.

\subsection{Elementary operations}\label{subsec23}

Consider a finite simply connected sub-graph $G$ of the
equilateral triangular lattice $\TT$, and let $\partial G$ be the
cycle of $G$ consisting of its boundary edges. Denote by
$\Q(\partial G)$ the set of triangular quadri-tilings whose
underlying tilings are lozenge tilings of $G$. Let $\Ls^1$ be a
lozenge tiling of $G$, and $L^1$ be the corresponding
lozenge-with-diagonals tiling. Then using the bijection between
the first height function and $2$-tilings of $L^1$ we obtain, in
exactly the same way as Elkies, Kuperberg, Larsen, Propp
\cite{EKLP} have for domino tilings, the following lemma:

\begin{lem}\label{lem4}
Every $2$-tiling of $L^1$ can be transformed into any other by a
finite sequence of the following operations, (in brackets is the
number of possible orientations for the graph corresponding to the
operation):
\begin{figure}[h]
\begin{center}
\includegraphics[width=\linewidth]{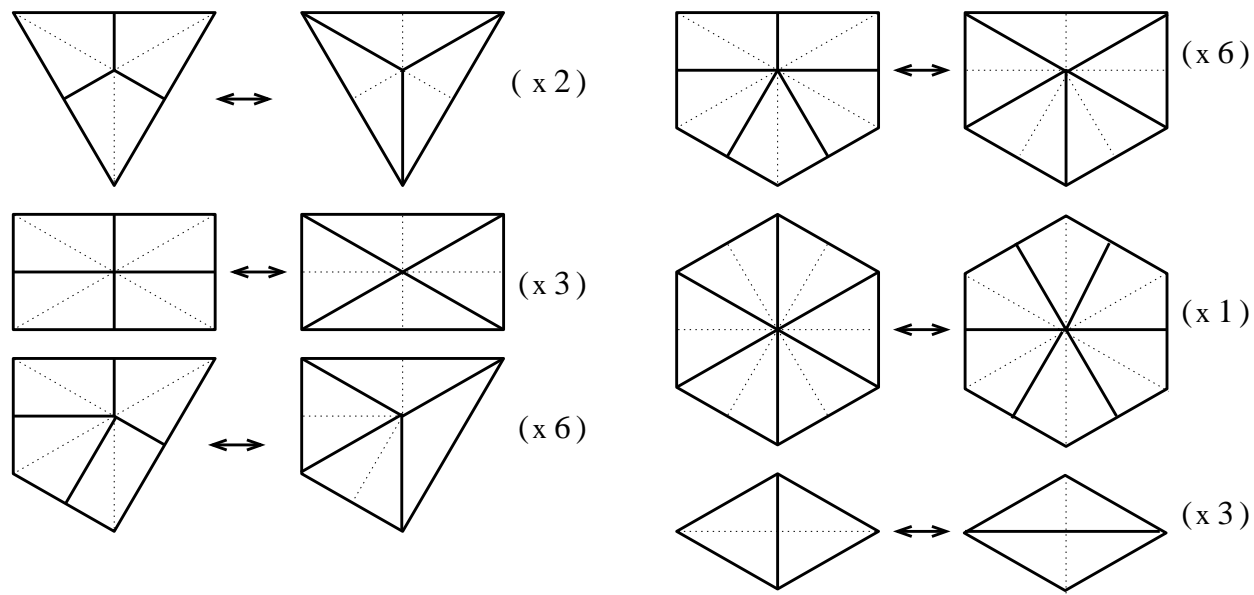}
\end{center}
\caption{Quadri-tile operations.} \label{fig4}
\end{figure}
 \end{lem}

\noindent Let us call {\bf quadri-tile operations} the 21
operations described in Lemma~\ref{lem4}. Moreover, every lozenge
tiling of $G$ can be transformed into any other by a finite
sequence of {\bf lozenge operations}:
\begin{center}
\includegraphics[height=2.3cm]{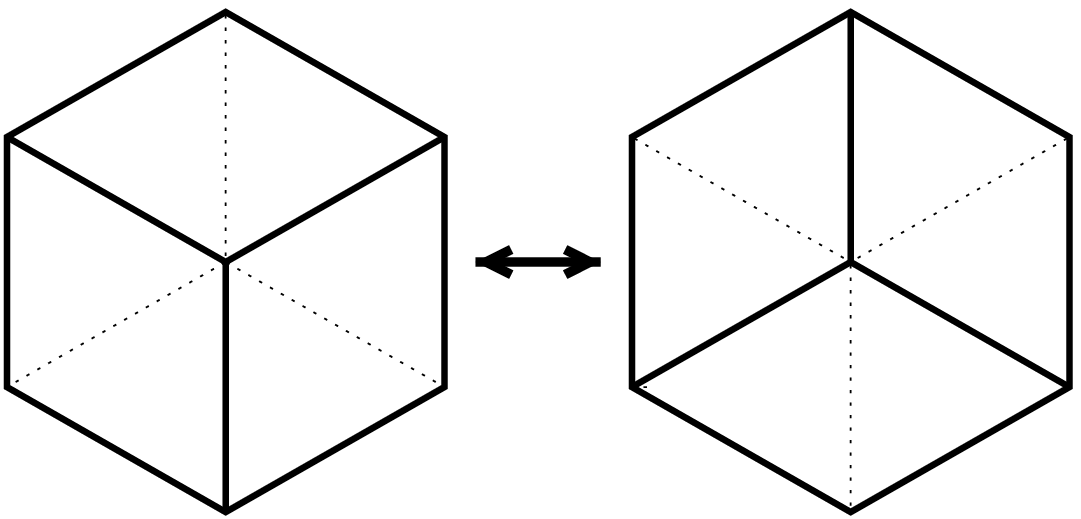}
\end{center}
\begin{figure}[h]
\caption{Lozenge operations.} \label{fig5}
\end{figure}
Note that if $\Ls$ is any lozenge tiling of $G$, then $L$ is
quadri-tilable with quadri-tiles obtained by cutting in two every
lozenge along one of its diagonals. Moreover, when one performs a
lozenge operation on such a quadri-tiling, one still obtains a
quadri-tiling of $\Q(\partial G)$. Let us call {\bf elementary
operations} the quadri-tile operations and the lozenge operations
performed on quadri-tilings as described above. Then we have:

\begin{lem}\label{lem5}
Every quadri-tiling of $\Q(\partial G)$ can be transformed into
any other by a finite sequence of elementary operations.
\end{lem}
\begin{proof}
This results from Lemma \ref{lem4}, and the above observation.
\end{proof}

\section{Gibbs measure on quadri-tilings}\label{sec3}

This section aims at giving a precise statement of Theorem
\ref{thm0} of the introduction (see Theorem \ref{thm2} below).
Sections \ref{subsec31} and \ref{subsec32} are taken form
\cite{Kenyon1} and give a precise definition of an isoradial
graph, the critical weight function, the Dirac and inverse Dirac
operator. Section \ref{subsec33} consists in the statement of
Theorem \ref{thm2}.

\subsection{Isoradial graphs and critical weight function}\label{subsec31}

\noindent The definition of the critical weight function follows
\cite{Kenyon1}. It is defined on edges of graphs satisfying a
geometric condition called {\bf isoradiality}: all faces of an
isoradial graph are inscribable in a circle, and all circumcircles
have the same radius.\\
Note that if $\Rs$ is a rhombus tiling of the plane, then the
corresponding rhombus-with-diagonals tiling $R$ is an isoradial
graph. Let us consider the embedding of the dual graph ${R}^*$
(the same notation is used for the one-skeleton of a graph and its
embedding) where the dual vertices are the circumcenters of the
corresponding faces. Then ${R}^*$ is also an isoradial graph and
the circumcenters of the faces are the
vertices of $R$.\\
To each edge $e$ of ${R}^*$, we associate a unit side-length
rhombus $R(e)$ whose vertices are the vertices of $e$ and the
vertices of its dual edge. Let $\widetilde{R}=\cup_{e \in {R}^*}
R(e)$. Note that the dual edges corresponding to the boundary
edges of the rhombi of $R$ have length zero, and that the rhombi
associated to these edges are degenerated.\\
\noindent For each edge $e$ of ${R}^*$, define $\nu(e)=2 \sin
\theta$, where $2\theta$ is the angle of the rhombus $R(e)$ at the
vertex it has in common with $e$; $\theta$ is called the {\bf
rhombus angle} of the edge $e$. Note that $\nu(e)$ is the length
of $e^*$, the dual edge of $e$. The function $\nu$ is called the
{\bf critical weight function}.

\subsection{Dirac and inverse Dirac operator}\label{subsec32}

\noindent Results and definitions of this section are due to
Kenyon \cite{Kenyon1}, see also Mercat \cite{Me2}. Note that this
section (as the previous one) is true for general isoradial graphs
with bipartite dual graphs.\jump

\noindent Let $\Rs$ be a rhombus tiling of the plane, then by
Lemma \ref{lem2}, ${R}^*$ is a bipartite graph. Let $B$ (resp.
$W$) be the set of black (resp. white) vertices of ${R}^*$, and
denote by $\nu$ the critical weight function on the edges of
${R}^*$. The Hermitian matrix $K$ indexed by the vertices of
${R}^*$ is defined as follows. If $v_1$ and $v_2$ are not adjacent
$K(v_1,v_2)=0$. If $w \in W$ and $b \in B$ are adjacent vertices,
then $K(w,b)=\overline{K(b,w)}$ is the complex number of modulus
$\nu(w b)$ and direction pointing from $w$ to $b$. If $w$ and $b$
have the same image in the plane, then $|K(w,b)|=2$, and the
direction of $K(w,b)$ is that which is perpendicular to the
corresponding dual edge, and has sign determined by the local
orientation. The infinite matrix $K$ defines the {\bf Dirac
operator} $K : \C^{V({R}^*)} \rightarrow \C^{V({R}^*)}$, by
$$
(Kf)(v)=\sum_{u \in {R}^*} K(v,u)f(u).
$$
where $V({R}^*)$ denotes the set of vertices of the graph
$R^*$.\jump

\noindent The {\bf inverse Dirac operator} $K^{-1}$ is defined
to be the operator satisfying:
\begin{enumerate}
\item $K K^{-1}=\mbox{Id}$,
\item $K^{-1}(b,w) \rightarrow 0$, when $|b-w| \rightarrow \infty$.
\end{enumerate}

\noindent Kenyon \cite{Kenyon1} obtains an explicit expression for
$K^{-1}$. Before stating his theorem, we need to define the
rational functions $f_{wv}(z)$. Let $w$ be a white vertex of
${R}^*$. For every other vertex $v$, define $f_{wv}(z)$ as
follows. Let $w=v_0,v_1,v_2,\ldots,v_k=v$ be an edge-path of
$\widetilde{R}$, from $w$ to $v$. Each edge $v_j v_{j+1}$ has
exactly one vertex of ${R}^*$ (the other is a vertex of $R$).
Direct the edge away from this vertex if it is white, and towards
this vertex if it is black. Let $e^{i\alpha_j}$ be the
corresponding vector in $\widetilde{R}$ (which may point contrary
to the direction of the path), then $f_{wv}$ is defined
inductively along the path, starting from
$$f_{ww}(z)=1.$$
If the edge leads away from a white vertex, or towards a black
vertex, then
$$
f_{w v_{j+1}}(z)=\frac{f_{w v_j}(z)}{z-e^{i\alpha_j}},
$$
else, if it leads towards a white vertex, or away from a black
vertex, then
$$
f_{w v_{j+1}}(z)=f_{w v_j}(z).(z-e^{i\alpha_j}).
$$

\noindent The function $f_{wv}(z)$ is well defined (i.e.
independent of the edge-path of $\widetilde{R}$ from $w$ to $v$),
because the multipliers for a path around a rhombus of
$\widetilde{R}$ come out to $1$. For a black vertex $b$ the value
$K^{-1}(b,w)$ will be the sum over the poles of $f_{wb}(z)$ of the
residue of $f_{wb}$ times the angle of $z$ at the pole. However,
there is an ambiguity in the choice of angle, which is only
defined up to a multiple of $2 \pi$. To make this definition
precise, angles are assigned to the poles of $f_{wb}(z)$. Working
on the branched cover of the plane, branched over $w$, so that for
each black vertex $b$ in this cover, a real angle $\theta_0$ is
assigned to the complex vector $b-w$, which increases by $2 \pi$
when $b$ winds once around $w$. In the branched cover of the
plane, a real angle in $[\theta_0-\pi+\Delta,\theta_0+\pi-\Delta]$
can be assigned to each pole of $f_{wb}$, for some small
$\Delta>0$.

\begin{thm}{\bf \cite{Kenyon1}}\label{thm3}
There exists a unique $K^{-1}$ satisfying the above two
properties, and $K^{-1}$ is given by:
$$
K^{-1}(b,w)=\frac{1}{4\pi^2 i}\int_{C} f_{wb}(z) \log z \;dz,
$$
where $C$ is a closed contour surrounding cclw the part of the
circle $\{ e^{i \theta}|\theta \in [\theta_0-\pi+
\Delta,\theta_0+\pi-\Delta]\}$, which contains all the poles of
$f_{wb}$, and with the origin in its exterior, see Figure \ref{fig16}.
\begin{figure}[h]
\begin{center}
\includegraphics[height=3.8cm]{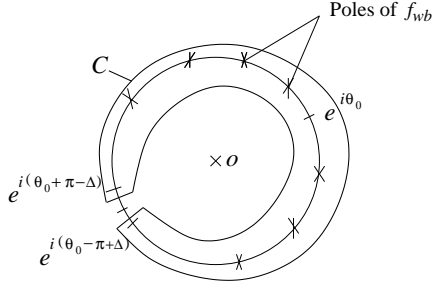}
\end{center}
\caption{An example of contour $C$.} \label{fig16}
\end{figure}
\end{thm}

\noindent The remarkable property of $K^{-1}(b,w)$ is that it only depends
on the local geometry of the graph, i.e. on an edge-path from $w$
to $b$.

\subsection{Statement of result}\label{subsec33}

\noindent Let $\Rs$ be a rhombus tiling of the plane, and $R$ be
the corresponding rhombus-with-diagonals tiling. Suppose that the
critical weight function $\nu$ is assigned to edges of $R^*$, and
let $K$ be the Dirac operator indexed by the vertices of $R^*$.
Moreover, recall that if $e_1=w_1 b_1,\ldots, e_k=w_k b_k$ is a
subset of edges of $R^*$, then the {\bf cylinder}
$\{e_1,\ldots,e_k\}$ is defined to be the set of dimer
configurations of $R^*$ which contain these edges. Let $\A$ be the
field consisting of the empty set and of the finite disjoint
unions of cylinders. Denote by $\sigma(\A)$ the $\sigma$-field
generated by $\A$.

\begin{thm}\label{thm2}
There is a probability measure $\mu^{R}$ on
$(\M({R}^*),\sigma(\A))$ such that for every cylinder
$\{e_1,\ldots,e_k\}$ of $R^*$,
\begin{equation}\label{GibbsM}
\mu^R(e_1,\ldots,e_k)=\left(\prod_{i=1}^{k} K(w_i,b_i)\right)
\det_{1 \leq i,\;j \leq k} \left(K^{-1}(b_i,w_j)\right).
\end{equation}
Moreover $\mu^R$ is a Gibbs measure on $\M(R^*)$. When $R^*$ is
doubly periodic, $\mu^R$ is the unique Gibbs measure which has
minimal free energy per fundamental domain among the two-parameter
family of ergodic Gibbs measures of \cite{KeOS}.
\end{thm}
From the proof, it appears that Theorem \ref{thm2} is in fact true
for all doubly periodic isoradial graphs with bipartite dual
graphs. We refer to the introduction for comments on Theorem \ref{thm2}.\vspace{0.6cm}

\section{Proof of Theorem \ref{thm2}}\label{sec4}

The proof of Theorem \ref{thm2} uses Propositions \ref{prop1} and
\ref{prop2} below. Proposition \ref{prop1} is a geometric property
of rhombus tilings and is the subject of Section \ref{subsec41}.
Proposition \ref{prop2} concerns the convergence of the Boltzmann
measure on some appropriate toroidal graphs, it is the subject of
Section \ref{subsec43}. In Section \ref{subsec42}, we introduce
the {\em real} Dirac operator and its inverse. This operator is
related to the Dirac operator and is needed for the proof of
Proposition \ref{prop2}. Theorem \ref{thm2} is then proved in
Section \ref{subsec44}.\vspace{1.5cm}

\subsection{Geometric property of rhombus tilings}\label{subsec41}

\begin{proposition}\label{prop1}
Let $\Rs$ be a rhombus tiling of the plane, then any finite simply
connected sub-graph $\Ps$ of $\Rs$ can be embedded in a periodic
rhombus tiling $\Ss$ of the plane.
\end{proposition}

\begin{proof}
This proposition is a direct consequence of Lemmas \ref{lem6},
\ref{lem7}, and Theorem \ref{lem8} below.
\end{proof}

\noindent The notion of {\bf train-track} has been introduced by
Mercat \cite{Me1}, see also Kenyon and Schlenker
\cite{Kenyon1,KeS}. A train-track of a rhombus tiling is a path of
rhombi (each rhombus being adjacent along an edge to the previous
rhombus) which does not turn: on entering a rhombus, it exits
across the opposite edge. Train-tracks are assumed to be maximal
in the sense that they extend in both directions as far as
possible. Thus train-tracks of rhombus tilings of the plane are
bi-infinite. Each rhombus in a train-track has an edge parallel to
a fixed unit vector $e$, called the {\bf transversal direction} of
the train-track. Let us denote by $t_e$ the train-track of
transversal direction $e$. In an oriented train-track (i.e. the
edges of the two parallel boundary paths of the train-track have
the same given orientation), we choose the direction of $e$ so
that when the train-track runs in the direction given by the
orientation, $e$ points from the right to the left. The vector $e$
is called the {\bf oriented transversal direction} of the oriented
train-track. A train-track cannot cross itself, and two different
train-tracks can cross at most once. A finite simply connected
sub-graph $\Ps$ of a rhombus tiling of the plane $\Rs$ is {\bf
train-track-convex}, if every train-track of $\Rs$ that intersects
$\Ps$ crosses the boundary of $\Ps$ twice exactly.\vspace{0.5cm}

\begin{lem}\label{lem6}
Let $\Rs$ be a rhombus tiling of the plane, then any finite simply
connected sub-graph $\Ps$ of $\Rs$ can be completed by a finite
number of rhombi of $\Rs$ in order to become train-track-convex.
\end{lem}

\begin{proof}
Let $e_1,\ldots,e_m$ be the boundary edges of $\Ps$. Every rhombus
of $\Ps$ belongs to two train-tracks of $\Rs$, each of which can
be continued in both directions up to the boundary of $\Ps$. In
both directions the intersection of each of the train-tracks and
the boundary of $\Ps$ is an edge parallel to the transversal
direction of the train-track. Thus, to take into account all
train-tracks of $\Rs$ that intersect $\Ps$, it suffices to
consider for every $i$ the train-tracks $t_{e_i}$ associated to
the boundary edges of $\Ps$. Consider the following algorithm (see
Figure \ref{fig6}).
\begin{quote}
Set $\Qs_1=\Ps$.\vspace{0.1cm}\\
For $i=1,\ldots,m$, do the following:\vspace{0.2cm}\\
Consider the train-track $t_{e_i}$, and let $2n_i$ be the number
of times $t_{e_i}$ intersects the boundary of $\Qs_i$.\vspace{0.2cm}
    \begin{itemize}
\item If $n_i > 1$: there are $n_i-1$ portions of $t_{e_i}$ that
are outside of $\Qs_i$, denote them by
${t_{e_i}}^1,\ldots,{t_{e_i}}^{n_i-1}$. Then, since $\Qs_i$ is
simply connected, for every $j$, $\Rs\, \backslash \, (\Qs_i \cup
{t_{e_i}}^j)$ is made of two disjoint sub-graphs of $\Rs$, one of
which is finite (it might be empty in the case where one of the
two parallel boundary paths of ${t_{e_i}}^j$ is part of the
boundary path of $\Qs_i$). Denote by ${g_{e_i}}^j$ the simply
connected sub-graph of $\Rs$ made of the finite sub-graph of
$\Rs\, \backslash \, (\Qs_i \cup {t_{e_i}}^j)$ and of
${t_{e_i}}^j$. Denote by ${b_{e_i}}^j$ the portion of the boundary
of $\Qs_i$ which bounds ${g_{e_i}}^j$. Replace $\Qs_i$ by
$\Qs_{i+1}= \Qs_i \cup (\cup_{j=1}^{n_i-1} {g_{e_i}}^j)$. By this
construction $t_{e_i}$ intersects the boundary of $\Qs_{i+1}$
exactly twice, and $\Qs_{i+1}$ is simply connected. \vspace{0.2cm}
\item If $n_{i}=1$: set $\Qs_{i+1}=\Qs_i$.\vspace{1cm}
    \end{itemize}
\end{quote}

\begin{figure}[h]
\begin{center}
\includegraphics[height=4.5cm]{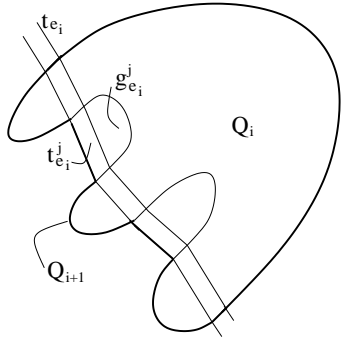}
\end{center}
\caption{One step of the algorithm.} \label{fig6}
\end{figure}

\noindent Let us show that at every step the train-tracks of $\Rs$
that intersect $\Qs_i$ and $\Qs_{i+1}$ are the same. By
construction, boundary edges of $\Qs_{i+1}$ are boundary edges of
$\Qs_i$ and of ${t_{e_i}}^j$, for every $j$. Let $f$ be an edge on
the boundary of $\Qs_{i+1}$, but not of $\Qs_i$, that is $f$ is on
the boundary of ${t_{e_i}}^j$ for some $j$, thus $t_f$ crosses
${g_{e_i}}^j$. Since two train-tracks cross at most once, $t_f$
has to intersect ${b_{e_i}}^j$, which means $t_f$ also crosses
$\Qs_i$. From this we also conclude that if a train-track
intersects the boundary of $\Qs_i$ twice, then it also intersects
the boundary of $\Qs_{i+1}$ twice.

\noindent Thus all train-tracks that intersect $\Qs_{m+1}$ cross
its boundary exactly twice, and $\Qs_{m+1}$ contains $\Ps$.
\end{proof}

\begin{lem}\label{lem7}
Let $\Rs$ be a rhombus tiling of the plane. Then any finite simply
connected train-track-convex sub-graph $\Ps$ of $\Rs$ can be
completed by a finite number of rhombi in order to become a convex
polygon $\Qs$, whose opposite boundary edges are parallel.
\end{lem}

\begin{proof}
Let $e_1,\ldots,e_m$ be the boundary edges of $\Ps$ oriented cclw.
Since $\Ps$ is train-track-convex, the train-tracks
$t_{e_1},\ldots,t_{e_m}$ intersect the boundary of $\Ps$ twice, so
that there are pairs of parallel boundary edges. Let us assume
that the transversal directions of the train-tracks are all
distinct (if this is not the case, one can always perturb the
graph a little so that it happens). Let us also denote by
$t_{e_1},\ldots,t_{e_m}$ the portions of the bi-infinite
train-tracks of $\Rs$ in $\Ps$. In what follows, indices will be
denoted cyclically, that is $e_j=e_{m+j}$. Write $x_j$ (resp.
$y_j$) for the initial (resp. end) vertex of an edge $e_j$.

\noindent Let $e_i, e_{i+1}$ be two adjacent boundary edges of
$\Ps$. Consider the translate $e_{i+1}^t$ of $e_{i+1}$ so that the
initial vertex of $e_{i+1}^t$ is adjacent to the initial vertex of
$e_i$. Then we define the {\bf turning angle from $e_i$ to
$e_{i+1}$} (also called exterior angle) to be the angle
$\widehat{e_i e_{i+1}^t}$, and we denote it by
$\theta_{e_i,e_{i+1}}$. If $e_i,e_j$ are two boundary edges, then
the {\bf turning angle from $e_i$ to $e_j$} is defined by
$\sum_{\alpha=i}^{j-1}\theta_{e_{\alpha},e_{\alpha+1}}$, and is
denoted by $\theta_{e_i,e_j}$.\vspace{0.5cm}\\

\noindent {\bf Properties}
\begin{enumerate}
\item $\sum_{\alpha=1}^{m}\theta_{e_{\alpha},e_{\alpha+1}}=2 \pi$.
\item If $e_i,e_j$ are two boundary edges, and if
$\gamma=\{f_1,\ldots,f_n\}$ is an oriented edge-path in $\Ps$ from
$y_i$ to $x_j$, then $\theta_{e_i,e_j}= \theta_{e_i,f_1} +
\sum_{\alpha=1}^{n-1}\theta_{f_{\alpha},f_{\alpha+1}}
+\theta_{f_n,e_j}$. \item If $e_i$ is a boundary edge of $\Ps$,
and $e_k$ is the second boundary edge at which $t_{e_i}$
intersects the boundary of $\Ps$, then $\theta_{e_i,e_k}=\pi$.
Thus $e_k$ and $e_i$ are oriented in the opposite direction, and
we denote $e_k$ by ${e_i}^{-1}$. \item $\Ps$ is convex, if and
only if every train-track of $\Ps$ crosses every other train-track
of $\Ps$.
\end{enumerate}

\noindent We first end the proof of Lemma \ref{lem7}, and then
prove Properties 1 to 4.\\
Note that Properties 1 and 2 are true for any finite simply
connected sub-graph of $\Rs$.\\ The number of train-tracks
intersecting $\Ps$ is $n=m/2$. So that if every train-track
crosses every other train-track, the total number of crossings is
$n(n-1)/2$. Consider the following algorithm (see Figure
\ref{fig7} for an example).
\begin{quote}
Set $\Qs_1=\Ps$, $n_1=$ the number of train-tracks that cross in
$\Qs_1$.\\ For $i=1,2,\ldots$ do the following:
\begin{itemize}
\item If $n_i=n(n-1)/2$: then by Property 4, $\Qs_i$ is convex.
\item If $n_i<n(n-1)/2$: then by Property 4,
$\theta_{e_{j_i},e_{j_i+1}}<0$ for some $j_i \in \{ 1,\ldots,m\}$.
Add the rhombus $\ell_{j_i}$ of parallel directions $e_{j_i},
e_{j_i+1}$ along the boundary of $\Qs_i$. Set $\Qs_{i+1}=\Qs_i
\cup \ell_{j_i}$, and rename the boundary edges $e_1,\ldots,e_m$
in cclw order. Then the number of train-tracks that cross in
$\Qs_{i+1}$ is $n_i+1$, set $n_{i+1}=n_i+1$. Note that if Property
4 is true for $\Qs_i$, it stays true for $\Qs_{i+1}$, and note
that the same train-tracks intersect $\Qs_i$ and $\Qs_{i+1}$.
\end{itemize}
\end{quote}\ \vspace{-1cm}\\
\begin{figure}[h]
\begin{center}
\includegraphics[height=5.7cm]{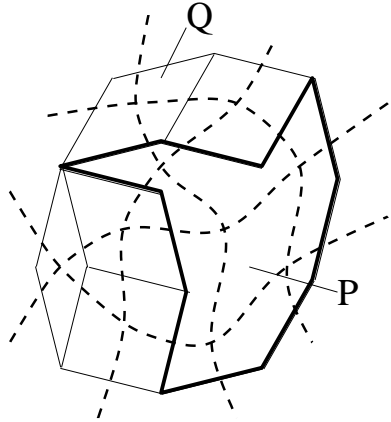}
\end{center}
\caption{Example of application of the algorithm.} \label{fig7}
\end{figure}

\noindent For the algorithm to be able to add the rhombus
$\ell_{j_i}$ at every step, we need to check that:
\begin{equation}\label{angle1}
\forall\; i,\;\theta_{e_{j_i-1},e_{j_i+1}}>-\pi, \mbox{ and
}\theta_{e_{j_i},e_{j_i+2}}>-\pi.
\end{equation}
\noindent Assume we have proved that for any finite simply
connected train-track-convex sub-graph $\Ps$ of $\Rs$ we have:
\begin{equation}\label{angle2}
\forall \; i,j,\; \theta_{e_i,e_j}>-\pi.
\end{equation}
\noindent Then Properties 1 and 2 imply that if (\ref{angle2})
is true for $\Qs_i$ it stays true for $\Qs_{i+1}$, moreover
(\ref{angle2}) implies (\ref{angle1}). So let us prove
(\ref{angle2}) by induction on the number of rhombi contained in
$\Ps$. If $\Ps$ is a rhombus, then (\ref{angle2}) is clear. Now
assume $\Ps$ is made of $k$ rhombi. Consider the train-tracks in
$\Ps$ adjacent to the boundary (every boundary edge $e$ of $\Ps$
belongs to a rhombus of $\Ps$ which has parallel directions $e$
and $f$; for every boundary edge $e$, the train-track of
transversal direction $f$ is the train-track adjacent to the
boundary). Denote the train-tracks adjacent to the boundary by
$t_1,\ldots,t_p$ in cclw order, and write $f_{\beta}$ for the
oriented transversal direction of $t_{\beta}$ (when the boundary
edge-path of $\Ps$ is oriented cclw). Consider two adjacent
boundary edges $e_i,e_{i+1}$ of $\Ps$ that don't belong to the
same boundary train-track. That is $e_i$ belongs to $t_{\beta}$,
and $e_{i+1}$ to $t_{\beta+1}$. Then either $\widehat{f_{\beta}
f_{\beta+1}}<0$ or $\widehat{f_{\beta} f_{\beta+1}}>0$, in the
second case $t_{\beta}$ and $t_{\beta+1}$ cross and their
intersection is a rhombus $\ell_{\beta}$ of $\Ps$. The rhombus
$\ell_{\beta}$ has boundary edges $e_i, e_{i+1}$, and
$f_{\beta+1}={e_i}^{-1},{f_{\beta}}^{-1}={e_{i+1}}^{-1}$. Now
Property 1 implies that $\sum_{\beta=1}^{p-1} \widehat{f_{\beta}
f_{\beta+1}}=2 \pi$, so that there always exists $\beta_0$ such
that $\widehat{f_{\beta_0} f_{\beta_0+1}} >0$. Removing
$\ell_{\beta_0}$ from $\Ps$ and using the assumption that $\Ps$ is
train-track-convex, we obtain a graph $\Ps'$ made of $k-1$ rhombi
which is train-track-convex. By induction, $\theta_{e,f}>-\pi$ for
every boundary edge of $\Ps'$, and using Property 2, we conclude
that this stays true for $\Ps$.\jump

\noindent Denote by $\Qs$ the convex polygon obtained from $\Ps$
by the algorithm, and assume that opposite boundary edges are not
parallel. Then there are indices $i$ and $j$ such that $e_i$ comes
before $e_j$, and $e_j^{-1}$ comes before $e_i^{-1}$. This implies
that $\theta_{e_i,e_j}=-\theta_{{e_j}^{-1},{e_i}^{-1}}$, so that
one of the two angles is negative, which means $\Qs$ can not be
convex. Thus we have a contradiction, and we conclude that
opposite boundary edges of $\Qs$ are parallel.\vspace{-0.2cm}\\

\noindent {\it Proof of Properties 1 to 4.}
\begin{description}
    \item 1. and 2. are straightforward.
    \item 3. When computing
$\theta_{e_i,e_k}$ along the boundary edge-path of $t_{e_i}$ we
obtain $\pi$, so by Property 2 we deduce that
$\theta_{e_i,e_k}=\pi$ in $\Ps$.
    \item 4. $\Ps$ is convex if and
only if, for every $i$, $\theta_{e_i,e_{i+1}}>0$, which is
equivalent to saying that, for every $i \neq j$,
$\theta_{e_i,e_j}>0$. Therefore Property 4 is equivalent to
proving that $\theta_{e_i,e_j}>0$, for every $i \neq j$, if and
only if every train-track of $\Ps$ crosses every
other train-track of $\Ps$.\\
\noindent Assume there are two distinct train-tracks $t_{e_\ell}$
and $t_{e_k}$ that don't cross in $\Ps$. Then, in cclw order
around the boundary of $\Ps$, we have either
$e_\ell,{e_k}^{-1},e_k,{e_\ell}^{-1}$, or
$e_\ell,e_k,{e_k}^{-1},{e_\ell}^{-1}$. It suffices to solve the
second case, the first case being similar. By Property 1,
$\theta_{e_\ell,e_k}+\theta_{e_k,{e_k}^{-1}}+
\theta_{{e_k}^{-1},{e_\ell}^{-1}}+\theta_{{e_\ell}^{-1},e_\ell}=2
\pi$. Moreover by Property 3,
$\theta_{e_k,{e_k}^{-1}}=\theta_{{e_\ell}^{-1},e_\ell}=\pi$, which
implies $\theta_{e_\ell,e_k}=-\theta_{{e_k}^{-1},{e_\ell}^{-1}}$.
Since all train-tracks have different transversal directions,
either $\theta_{e_\ell,e_k}$ or
$\theta_{{e_k}^{-1},{e_\ell}^{-1}}$ is negative.\\
Now take two
boundary edges $e_i,e_j$ of $\Ps$ (with $\;i \neq j$, and $e_j
\neq {e_i}^{-1}$), and assume the train-tracks $t_{e_i}, t_{e_j}$
cross inside $\Ps$. Then in cclw order around the boundary of
$\Ps$, we have either $e_i,{e_j}^{-1},{e_i}^{-1},e_j$, or
$e_i,e_j,{e_i}^{-1},{e_j}^{-1}$. It suffices to solve the second
case since the first case can be deduced from the second one. The
intersection of $t_{e_i}$ and $t_{e_j}$ is a rhombus $\ell$. Let
${\tilde{e}_j}^{-1}$ (resp. ${\tilde{e}_i}^{-1}$) be the boundary
edge of $\ell$ parallel and closest to $e_j$ (resp. $e_i$),
oriented in the opposite direction, then
$\theta_{{\tilde{e}_j}^{-1},{\tilde{e}_i}^{-1}}<0$. Let $\gamma_j$
(resp. $\gamma_i$) be the boundary edge-path of $t_{e_j}$ (resp.
$t_{e_i}$) from $y_j$ to $\tilde{x}_j$ (resp. from $\tilde{y}_i$
to $x_i$), and let $\Qs$ be the sub-graph of $\Rs$ whose boundary
is
$e_i,e_{i+1},\ldots,e_j,\gamma_j,{\tilde{e}_j}^{-1},{\tilde{e}_i}^{-1},
\gamma_i$. Since $t_{e_i}$ and $t_{e_j}$ intersect the boundary of
$\Ps$ twice, they also intersect the boundary of $\Qs$ twice.
Moreover $t_{e_i}$ and $t_{e_j}$ don't cross in $\Qs$, so that
$\theta_{e_i,e_j}=-
\theta_{{\tilde{e}_j}^{-1},{\tilde{e}_i}^{-1}}>0$.
\hspace*{\fill} $\square$
\end{description}
\end{proof}

\begin{thm}\label{lem8}
Any convex $2n$-gon $\Qs$ whose opposite boundary edges are
parallel and of the same length can be embedded in a periodic
tiling of the plane by $\Qs$ and rhombi.
\end{thm}

\begin{proof}
\noindent Let $e_1,\ldots,e_n,{e_1}^{-1},\ldots,{e_n}^{-1}$ be the
boundary edges of the polygon $\Qs$ oriented cclw. If $n \leq 3$,
then $\Qs$ is either a rhombus or a hexagon, and it is
straightforward that the plane can be tiled periodically with
$\Qs$.

\noindent If $n>4$, for $k=1,\ldots,n-3$, do the following (see
Figure \ref{fig8}): along $e_{n-k}$ add the finite train-track
$\tilde{t}_{e_{n-k}}$ of transversal direction $e_{n-k}$, going
away from $\Qs$, whose boundary edges starting from the boundary
of $\Qs$ are:
$$
e_1,\underbrace{e_2,e_1},\underbrace{e_3,e_2,e_1},\ldots, \underbrace{e_{n-k-2},\ldots,e_1}.
$$ \ \vspace{0cm}\\

\begin{figure}[h]
\begin{center}
\includegraphics[height=5.1cm]{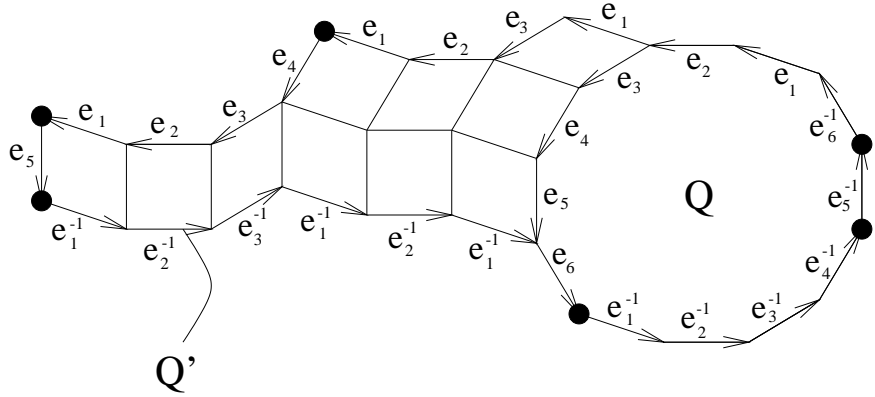}
\end{center}
\caption{Fundamental domain of a periodic tiling of the plane by
dodecagons and rhombi.} \label{fig8}
\end{figure}

\noindent Since the polygon $\Qs$ is convex, the rhombi that are
added are well defined, moreover the intersection of
$\tilde{t}_{e_i}$ and the boundary of $\Qs$ is the edge $e_i$, and
$\tilde{t}_{e_i}$ doesn't cross $\tilde{t}_{e_j}$ when $i \neq j$.
So we obtain a new polygon $\Qs'$ made of $\Qs$ and rhombi, whose
boundary edge-path is $\gamma_1,\ldots,\gamma_6$, (when starting
from the edge ${e_n}^{-1}$ of $\Qs$), where:
\begin{description}
\item $\gamma_1={e_n}^{-1},e_1,\underbrace{e_2,e_1}, \underbrace{e_3,e_2,e_1},\ldots,\underbrace{e_{n-3},\ldots,e_1},$
\item $\gamma_2=e_{n-2},\ldots,e_1,$
\item $\gamma_3=e_{n-1},$
\item $\gamma_4=\underbrace{{e_1}^{-1},\ldots,{e_{n-3}}^{-1}},\ldots,\underbrace{{e_1}^{-1},{e_2}^{-1}},{e_1}^{-1},e_n,$
\item $\gamma_5={e_1}^{-1},\ldots,{e_{n-2}}^{-1},$
\item $\gamma_6={e_{n-1}}^{-1}.$
\end{description}
Noting that
$\gamma_4={\gamma_1}^{-1},\gamma_5={\gamma_2}^{-1},\gamma_6={\gamma_3}^{-1}$,
and using the fact that the plane can be tiled periodically with
hexagons which have parallel opposite boundary edges, we deduce
that the plane can be tiled with $\Qs'$, that is it can be tiled
periodically by $\Qs$ and rhombi.
\end{proof}

\begin{remark}
It was pointed out by the referee that Theorem \ref{lem8} might be
known already. After a second look at the extensive literature on
tilings, we were not able to find a reference, but any information is
of course welcome.
\end{remark}

\subsection{Real Dirac and inverse real Dirac operator}\label{subsec42}

\noindent In the whole of this section we let $\Rs$ be a rhombus
tiling of the plane, and $R$ be the corresponding
rhombus-with-diagonals tiling. Assume that the critical weight
function $\nu$ is assigned to edges of $R^*$, and denote
by $K$ the Dirac operator indexed by the vertices of $R^*$.\\
The proof of Theorem \ref{thm2} requires to take the limit of
Boltzmann measures on some appropriate toroidal graphs (see
Section \ref{subsec43}). In order to do this, we need to introduce
the {\bf real Dirac operator} denoted by $\Ks$. Both the Dirac
operator $K$ and the real Dirac operator $\Ks$ are represented by
infinite weighted adjacency matrices indexed by the vertices of
$R^*$. For $K$, the edges of $R^*$ are un-oriented and weighted by
their critical weight times a complex number of modulus $1$. For
$\Ks$, edges of $R^*$ are oriented with a {\em clockwise odd}
orientation, and are weighted by their critical weight. Both
weight functions yield the same probability measure on finite
simply connected sub-graphs of $R^*$, but, and this is the reason
why we introduce the real Dirac operator, these weights {\em do
not} yield the same probability distribution on toroidal
sub-graphs of $R^*$.\\
The structure of this section is close to that of Section
\ref{subsec32}. We first define the real Dirac operator $\Ks$ and
its inverse $\Ks^{-1}$ and then, using results of \cite{Kenyon1},
we prove the existence and uniqueness of $\Ks^{-1}$ by giving an
explicit expression for $\Ks^{-1}$. Note that this section is
actually true for general isoradial graphs with bipartite dual
graphs.

\subsubsection{Real Dirac operator}\label{subsubsec421}\ \jump

\noindent Let us define an orientation of the edges of $R^*$. An
elementary cycle $C$ of $R^*$ is said to be {\bf clockwise odd}
if, when traveling cw around the edges of $C$, the number of
co-oriented edges is odd. Note that since $R^*$ is bipartite, the
number of contra-oriented edges is also odd. Kasteleyn
\cite{Kast2} defines the orientation of the graph $R^*$ to be {\bf
clockwise odd} if all elementary cycles are clockwise odd. He also
proves that, for planar simply connected graphs, such an
orientation always exists.\\
Consider a clockwise odd orientation of the edges of $R^*$. Define
$\Ks$ to be the infinite adjacency matrix of the graph $R^*$,
weighted by the critical weight function $\nu$. That is, if $v_1$
and $v_2$ are not adjacent, $\Ks(v_1,v_2)=0$. If $w\in W$ and
$b\in B$ are adjacent vertices, then
$\Ks(w,b)=-\Ks(b,w)=(-1)^{\II_{(w,b)}}\nu(wb)$, where
$\II_{(w,b)}=0$ if the edge $wb$ is oriented from $w$ to $b$, and
$1$ if it is oriented from $b$ to $w$. The infinite matrix $\Ks$
defines the {\bf real Dirac operator} $\Ks$: $\C^{V(R^*)}
\rightarrow \C^{V(R^*)}$, by
\begin{equation*}\label{618}
(\Ks f)(v)=\sum_{u \in R^*} \Ks(v,u)f(u),
\end{equation*}
The matrix $\Ks$ is also called a {\bf Kasteleyn matrix} for the
underlying dimer model.

\subsubsection{Inverse real Dirac operator}\label{subsubsec422}\ \jump

\noindent The {\bf inverse real Dirac operator} $\Ks^{-1}$ is
defined to be the unique operator satisfying
\begin{enumerate}
    \item $\Ks \Ks^{-1}=$ Id,
    \item $\Ks^{-1}(b,w)\rightarrow 0$, when $|b-w|\rightarrow\infty$.
\end{enumerate}
Let us define the rational functions $\fs_{wx}(z)$. They are the
analogous of the rational functions $f_{wv}(z)$, but are defined
for vertices $x\in R^*$ (whereas the functions $f_{wv}(z)$ were
defined for vertices $v\in \widetilde{R}$). Let $w\in W$, and let
$x\in B$ (resp. $x\in W$); consider the edge-path
$w=w_1,b_1,\ldots,w_k, b_k=x$ (resp.
$w=w_1,b_1,\ldots,w_k,b_k,w_{k+1}=x$) of $R^*$ from $w$ to $x$.
Let $R(w_j b_j)$ be the rhombus associated to the edge $w_j b_j$,
and denote by $w_j,x_j,b_j,y_j$ its vertices in cclw order;
$e^{i\alpha_j}$ is the complex vector $y_j-w_j$, and
$e^{i\beta_j}$ is the complex vector $x_j-w_j$. In a similar way,
denote by $w_{j+1},x_j',b_j,y_j'$ the vertices of the rhombus
$R(w_{j+1}b_j)$ in cclw order, then $e^{i\alpha_j'}$ is the
complex vector $y_j'-w_{j+1}$, and $e^{i\beta_j'}$ is the complex
vector $x_j'-w_{j+1}$. Then $\fs_{wx}(z)$ is defined inductively
along the path,
\begin{eqnarray*}\label{61}
\fs_{ww}(z)&=&1,\\
\fs_{w b_j}(z)&=&\fs_{w w_j}(z)
\frac{(-1)^{\II_{(w_j,b_j)}}e^{i\frac{\alpha_j+\beta_j}{2}}}{(z-e^{i\alpha_j})(z-e^{i\beta_j})},\\
\fs_{w w_{j+1}}(z)&=&\fs_{w b_j}
(z)(-1)^{\II_{(w_{j+1},b_j)}}e^{-i\frac{\alpha_j'+\beta_j'}{2}}(z-e^{i\alpha_j'})(z-e^{i\beta_j'}).\\
\end{eqnarray*}

\begin{remark}\label{rem1}\ \jump

\noindent We have the following relations between the real and the
complex case.
\begin{enumerate}
    \item $\forall\,w\in W,\;\forall\,x\in B\cup W$,
    $\fs_{wx}(z)=\overline{\fs_{wx}(0)}f_{wx}(z)$.
    \item $\forall\,w\in W,\;\forall\,b\in B,$ such that $w$ is
    adjacent to $b$,
    $\Ks(w,b)=\fs_{wb}(0)K(w,b)$.
\end{enumerate}
\end{remark}
\begin{proof}\ \jump

\noindent $1.$ This is a direct consequence of the definitions of the
functions $\fs_{wx}$ and $f_{wx}$.\jump

\noindent $2.$ Let $R(wb)$ be the rhombus associated to the edge $wb$,
and let $w,x,b,y$ be its vertices in cclw order. Denote by
$e^{i\alpha}$ the complex vector $y-w$, and by $e^{i\beta}$ the
complex vector $x-w$. Let $\theta$ be the rhombus angle of the
edge $wb$. By definition we have,
    \begin{eqnarray*}
    \Ks(w,b)&=&(-1)^{\II_{(w,b)}}2\sin\theta=(-1)^{\II_{(w,b)}}\frac{e^{i\frac{\alpha-\beta}{2}}-
    e^{-i\frac{\alpha-\beta}{2}}}{i},\\
    K(w,b)&=&i(e^{i\beta}-e^{i\alpha}),\\
    \fs_{wb}(0)&=&(-1)^{\II_{(w,b)}}e^{-i\frac{\alpha+\beta}{2}}.
    \end{eqnarray*}
Combining the above three equations yields $2$.
\end{proof}
\begin{lem}\label{lem9}
The function $\fs_{wx}$ is well defined.
\end{lem}
\begin{proof}
Showing that the function $\fs_{wx}$ is well defined amounts to
proving that $\fs_{wx}$ is independent of the edge-path of $R^*$
from $w$ to $x$. This is equivalent to proving the following: let
$w_1,b_1,\ldots,w_k,b_k,w_{k+1}=w_1$ be the vertices of an
elementary cycle $C$ of $R^*$, where vertices are enumerated in
cclw order; if $\fs_{w_1 w_1}(z)=1$ then $\fs_{w_1 w_{k+1}}(z)=1$.
Let us use the notations introduced in the definition of
$\fs_{wx}$, and denote indices cyclically, that is $k+1\equiv 1$.
By Remark \ref{rem1}, we have
\begin{equation*}
\fs_{w_1 w_{k+1}}(z)=\overline{\fs_{w_1 w_{k+1}}(0)}f_{w_1
w_{k+1}}(z).
\end{equation*}
Since the function $f_{wx}$ is well defined, $f_{w_1
w_{k+1}}(z)=1$. Hence, it remains to prove that
$\overline{\fs_{w_1 w_{k+1}}(0)}=1$. By definition of $\fs_{wx}$,
we have
\begin{equation*}\label{62}
\overline{\fs_{w_1 w_{k+1}}(0)}=\prod_{j=1}^k\left(
(-1)^{(\II_{(w_j,b_j)}+\II_{(w_{j+1},b_j)})}
e^{i\frac{\alpha_j+\beta_j}{2}}e^{-i\frac{\alpha_j'+\beta_j'}{2}}\right).
\end{equation*}
Moreover for every $j$, $\alpha_j'=\beta_j$, so that
\begin{equation*}\label{63}
\overline{\fs_{w_1 w_{k+1}}(0)}=\prod_{j=1}^k\left(
(-1)^{(\II_{(w_j,b_j)}+\II_{(w_{j+1},b_j)})}
e^{i\frac{\alpha_j-\beta_j}{2}}e^{i\frac{\alpha_j'-\beta_j'}{2}}\right).
\end{equation*}
Let $\theta_j$ (resp. $\theta_j'$) be the rhombus angle of the
edge $w_j b_j$ (resp. $w_{j+1}b_j$), then
\begin{equation}\label{64}
\overline{\fs_{w_1 w_{k+1}}(0)}=(-1)^{\sum_{j=1}^k
(\II_{(w_j,b_j)}+\II_{(w_{j+1},b_j)})} e^{i\sum_{j=1}^k
(\theta_j+\theta_j')}.
\end{equation}
The cycle $C$ corresponds to a face of the graph $R^*$. Let $c$ be
the circumcenter of this face, and let $\tau_j$ (resp. $\tau_j'$)
be the angle of the rhombus $R(w_j b_j)$ (resp. $R(w_{j+1}b_j)$)
at the vertex $c$. Then $\tau_j=\pi-2\theta_j$, and
$\tau_j'=\pi-2\theta_j'$. Since $\sum_{j=1}^k
(\tau_j+\tau_j')=2\pi$, we deduce
$\sum_{j=1}^k(\theta_j+\theta_j')=\pi(k-1)$. Hence,
\begin{equation}\label{65}
e^{i\sum_{j=1}^k (\theta_j+\theta_j')}=-(-1)^k.
\end{equation}
Moreover $\II_{(w_{j+1},b_j)}=1-\II_{(b_j,w_{j+1})}$, and
$(-1)^{1-\II_{(b_j,w_{j+1})}}=(-1)^{\II_{(b_j,w_{j+1})}-1}$, so
\begin{equation*}\label{66}
(-1)^{\sum_{j=1}^k (\II_{(w_j,b_j)}+\II_{(w_{j+1},b_j)})}=
(-1)^{\sum_{j=1}^k (\II_{(w_j,b_j)}+\II_{(b_j,w_{j+1})})-k}.
\end{equation*}
Note that $\sum_{j=1}^k (\II_{(w_j,b_j)}+\II_{(b_j,w_{j+1})})$ is
the number of co-oriented edges encountered when traveling cclw
around the cycle $C$. Since the orientation of the edges of $R^*$
is clockwise odd, it is also counterclockwise odd, and so this
number is odd. This implies
\begin{equation}\label{67}
(-1)^{\sum_{j=1}^k
(\II_{(w_j,b_j)}+\II_{(w_{j+1},b_j)})}=-(-1)^{-k}.
\end{equation}
The proof is completed by combining equations (\ref{64}),
(\ref{65}) and (\ref{67}).
\end{proof}
As in the complex case, a real angle in
$[\theta_0-\pi+\Delta,\theta_0+\pi-\Delta]$ can be assigned to
each pole of $\fs_{wb}$, for some small $\Delta>0$; where
$\theta_0$ is the real angle assigned to the vector $b-w$.

\begin{lem}\label{lem13}
There exists a unique $\Ks^{-1}$ satisfying the above two
properties, and $\Ks^{-1}$ is given by:
\begin{equation}\label{68}
\Ks^{-1}(b,w)=\frac{1}{4\pi^2 i}\int_{C} \fs_{wb}(z) \log z \;dz,
\end{equation}
where $C$ is a closed contour surrounding cclw the part of the
circle $\{ e^{i \theta} \,|\, \theta \in [\theta_0 - \pi +
\Delta,\theta_0 + \pi - \Delta ]\}$, which contains all the poles
of $\fs_{wb}$, and with the origin in its exterior.
\end{lem}
\begin{proof}
Let $\Fs(b,w)$ be the right hand side of (\ref{68}). Fix a vertex
$w_0\in W$, and let us prove that $\sum_{b\in
B}\Ks(w_0,b)\Fs(b,w)=\delta_{w_0}(w)$. Denote by $b_1,\ldots,b_k$
the black neighbors of $w_0$. Using Remark \ref{rem1}, we obtain
for every $j$,
\begin{eqnarray*}
\Ks(w_0,b_j)=\fs_{w_0 b_j}(0)K(w_0,b_j),\\
\fs_{w b_j}(z)=\overline{\fs_{w b_j}(0)}f_{w b_j}(z).
\end{eqnarray*}
Moreover $\forall\,w\in W,\,\forall\,b\in B$, we have
$\overline{\fs_{wb}(0)}=\fs_{w b}(0)^{-1}=\fs_{bw}(0)$, so that
$\fs_{w_0 b_j}(0)\overline{\fs_{w b_j}(0)}=\fs_{w_0
b_j}(0)\fs_{b_j w}(0)=\fs_{w_0 w}(0)$. Hence, using Theorem
\ref{thm3}, we obtain for every $j$,
\begin{equation*}
\Ks(w_0,b_j)\Fs(b_j,w)=\fs_{w_0 w}(0)K(w_0,b_j)K^{-1}(b_j,w).
\end{equation*}
Since $\Ks(w_0,b)=0$ when $w_0$ and $b$ are not adjacent, and
since $K^{-1}$ is the inverse Dirac operator, we obtain
{\small\begin{eqnarray*}\label{610}
    \sum_{b\in B}\Ks(w_0,b)\Fs(b,w)&=&\sum_{j=1}^k \Ks(w_0,b_j)\Fs(b_j,w),\\
        &=&\fs_{w_0 w}(0)\sum_{j=1}^k K(w_0,b_j)K^{-1}(b_j,w)=\fs_{w_0
w}(0)\delta_{w_0}(w)=\delta_{w_0}(w).
\end{eqnarray*}}
Uniqueness of $\Ks^{-1}$ follows from the uniqueness of $K^{-1}$.
\end{proof}

\subsection{Convergence of the Boltzmann measure on the
torus}\label{subsec43}

\noindent Let $\Rs$ be a rhombus tiling of the plane, and $R$ be
the corresponding rhombus-with-diagonals tiling. Suppose that the
critical weight function is assigned to edges of $R^*$, and denote
by $K$ the Dirac operator indexed by the vertices of $R^*$.\\
Consider a subset of edges $e_1=w_1 b_1,\ldots,e_k=w_k b_k$ of
$R^*$, and let $\Ps$ be a finite simply connected sub-graph of
$\Rs$ such that $P^*$ contains these edges. By Proposition
\ref{prop1}, there exists a periodic rhombus tiling of the plane
$\Ss$ that contains $\Ps$. Let $S$ be the corresponding
rhombus-with-diagonals tiling, and assign the critical weight
function to edges of $S^*$. Denote by $\Lambda$ the lattice which
acts periodically on $S$, and suppose that the dual graph
$\Snbar^*$ of the toroidal graph $\Snbar=S/n\Lambda$ is bipartite
(this is possible by eventually replacing $\Lambda$ by
$2\Lambda$). Denote by $\mu_n^S$ be the Boltzmann measure on dimer
configurations $\M(\Snbar^*)$ of $\Snbar^*$. Then we have,\\

\begin{proposition}\label{prop2}
\begin{equation*}
\lim_{n\rightarrow\infty}\mu_n^S(e_1,\ldots,e_k)=
    \left(\prod_{i=1}^k K(w_i,b_i)\right)\det_{1\leq i,j\leq
    k}K^{-1}(b_i,w_j).
\end{equation*}
\end{proposition}

\begin{proof}
Let us first define an orientation of the edges of $S^*$, and the
four Kasteleyn matrices $\Ks_1^n,\ldots,\Ks_4^n$ of the graph
$\Snbar^*$. Consider the graph $\S1bar^*$, then it is a bipartite
graph on the torus. Fix a reference matching $M_0$ of $\S1bar^*$.
For every other perfect matching $M$ of $\S1bar^*$, consider the
superposition $M\cup M_0$ of $M$ and $M_0$, then $M\cup M_0$
consists of doubled edges and cycles. Let us define four parity
classes for perfect matchings $M$ of $\S1bar^*$: (e,e) consists of
perfect matchings $M$ for which cycles of $M\cup M_0$ circle the
torus an even number of times horizontally and vertically; (e,o)
consists of perfect matchings $M$, for which cycles of $M\cup M_0$
circle the torus an even number of times horizontally, and an odd
number of times vertically; (o,e) and (o,o) are defined in a
similar way. By Tesler \cite{Tesler}, one can construct an
orientation of the edges of $\S1bar^*$, so that the corresponding
weighted adjacency matrix $\Ks_1^1$ has the following property:
perfect matchings which belong to the same parity class have the
same sign in the expansion of the determinant of $\Ks_1^1$. By an
appropriate choice of sign, we can make the (e,e) class have the
plus sign in $\det\Ks_1^1$, and the other three have minus sign.
Consider a horizontal and a vertical cycle of $\S1bar$. Then
define $\Ks_2^1$ (resp. $\Ks_3^1$) to be the matrix $\Ks_1^1$
where the sign of the coefficients corresponding to edges crossing
the horizontal (resp. vertical) cycle is reversed; and define
$\Ks_4^1$ to be the matrix $\Ks_1^1$ where the sign of the
coefficients corresponding to
edges crossing both cycles are reversed.\\
The orientation of the edges of $\S1bar^*$ defines a periodic
orientation of the graph $S^*$. For every $n$, consider the graph
$\Snbar^*$ and the four matrices $\Ks_1^n,\Ks_2^n,\Ks_3^n,\Ks_4^n$
defined as above. These matrices are called the {\bf Kasteleyn
matrices} of the graph $\Snbar^*$.\\
The orientation defined on the edges of the graph $S^*$ is a
clockwise odd orientation. Let $\Ks_S$ be the real Dirac operator
indexed by the vertices of $S^*$ corresponding to this clockwise
orientation, and let $K_S$ be the Dirac operator indexed by the
vertices of $S^*$. Then Proposition \ref{prop2} is a direct
consequence of Lemmas \ref{lem10}, \ref{lem11}, \ref{lem12} below.

\begin{lem}\label{lem10}
\begin{equation}\label{6}
\lim_{n\rightarrow\infty} \mu_n^S(e_1,\ldots,e_k)=
    \left(\prod_{i=1}^k \Ks_S(w_i,b_i)\right)\det_{1\leq i,j\leq
    k}\Ks_S^{-1}(b_i,w_j).
\end{equation}
\end{lem}
\begin{proof}
The {\bf toroidal partition function} $Z(\Snbar^*,\nu)$ is defined
to be the weighted sum (weighted by the function $\nu$) of dimer
configurations of the graph $\Snbar^*$. Then, by Tesler
\cite{Tesler} (it is a generalization of a theorem of Kasteleyn
\cite{Kast1}), we have
\begin{thm}{\rm \cite{Kast1,Tesler}}\label{thm4}
\begin{equation*}\label{4}
Z(\Snbar^*,\nu)=\frac{1}{2}(-\det\Ks_1^n+\det\Ks_2^n+\det\Ks_3^n+\det\Ks_4^n).
\end{equation*}
\end{thm}
Kenyon gives the following theorem for the Boltzmann measure
$\mu_n^S$.
\begin{thm}{\rm \cite{Kenyon}}\label{thm5}
$\mu_n^S(e_1,\ldots,e_k)$ is equal to
$\left(\prod_{i=1}^{k}\Ks_S(w_i,b_i)\right)$ times
{\scriptsize\begin{equation}\label{5}
\left(-\frac{\det\Ks_1^n}{2Z(\Snbar^*,\nu)}
    \det_{1 \leq i,j \leq k}\left( (\Ks_1^n)^{-1}(b_i,w_j)\right)
    + \sum_{\ell=2}^{4}\frac{ \det\Ks_\ell^n}{2Z(\Snbar^*,\nu)}
    \det_{1 \leq i,j \leq k}\left((\Ks_\ell^n)^{-1}(b_i,w_j)\right)\right).
\end{equation}}
\end{thm}
This part of the argument can be found in \cite{Kenyon}. Equation
(\ref{5}) is a weighted average of the four quantities $\det_{1
\leq i,j \leq k}\left((\Ks_\ell^n)^{-1}(b_i,w_j)\right)$, with
weights\\ $\frac{1}{2}\det \Ks_\ell^n/Z(\Snbar^*,\nu)$. These
weights are all in the interval $(-1,1)$ since, for every
$\ell=1,\ldots,4$, $2Z(\Snbar^*,\nu)>|\det \Ks_\ell^n|$. Indeed,
$Z(\Snbar^*,\nu)$ counts the weighted sum of dimer configurations
of $\Snbar^*$, whereas $|\det \Ks_\ell^n|$ counts some
configurations with negative sign. Moreover, by Theorem \ref{thm4}
these weights sum to $1$, so that the weighted average converges
to the same value as each
{\small $\displaystyle \det_{1 \leq i,j \leq k}\left((\Ks_\ell^n)^{-1}(b_i,w_j)\right)$}.\jump

\noindent Denote by $B_S$ (resp. $W_S$) the set of black (resp.
white) vertices of $S^*$. Let us prove that for every
$\ell=1,\ldots,4$, and for every $w\in W_S$, $b\in B_S$,
$(\Ks_\ell^n)^{-1}(b,w)$ converges to $\Ks_S^{-1}(b,w)$ on a
subsequence of $n'$s. The following theorem of \cite{KeOS} gives
the convergence on a subsequence of $n'$s of the inverse Kasteleyn
matrices of the graph $\Snbar^*$.
\begin{thm}{\rm\cite{KeOS}}\label{thm6}
For every $w\in W_S, b\in B_S$, $\ell=1,\ldots,4$,
\begin{equation}\label{340}
{\lim}'_{\,n \rightarrow
\infty}(\Ks_{\ell}^n)^{-1}(b,w)=\frac{1}{(2 \pi)^2} \int_{S^1
\times S^1} \frac{Q_{b,w}(z,u)u^x z^y}{P(z,u)} \frac{dz}{z}
\frac{du}{u},
\end{equation}
where $Q_{b,w}$ and $P$ are polynomials ($Q_{b,w}$ only depends on
the equivalence class of $w$ and $b$), and $x$ (resp. $y$) is the
horizontal (resp. vertical) translation from the fundamental
domain of $b$ to the fundamental domain of $w$.
\end{thm}
Denote by $\Fs(b,w)$ the right hand side of (\ref{340}). In
\cite{KeOS}, it is proved that $\Fs(b,w)$ converges to $0$ as
$|b-w|\rightarrow\infty$, as long as the dimer model is not in its
frozen phase. Moreover, it is proved in \cite{KO} that dimer
models on isoradial graphs are never in their frozen phase when
the weight function is the critical one. Hence, we deduce that
$\Fs(b,w)$ converges to $0$ as $|b-w|\rightarrow\infty$.\\
Let us prove that for every $b\in B_S,w\in W_S$,
$\Fs(b,w)=\Ks_S^{-1}(b,w)$. Consider $w_1,w_2\in W_S$, and denote
by $b_1,\ldots,b_k$ the neighbors of $w_1$. Assume $n$ is large
enough so that the graph $\Snbar^*$ contains
$w_1,w_2,b_1,\ldots,b_k$, and so that the edges $w_1 b_j$ do not
cross the horizontal and vertical cycle of $\Snbar$. Then, for
every $\ell=1,\ldots,4,$
\begin{equation*}\label{341}
\sum_{b\in B_S}
\Ks_\ell^n(w_1,b)(\Ks_\ell^n)^{-1}(b,w_2)=\sum_{j=1}^k
\Ks_\ell^n(w_1,b_j)(\Ks_\ell^n)^{-1}(b_j,w_2)=\delta_{w_1 w_2}.
\end{equation*}
Moreover, $\Ks_\ell^n(w_1,b_j)=\Ks_S(w_1,b_j)$, so that taking the
limit on a subsequence of $n'$s, and using Theorem \ref{thm6}, we
obtain
\begin{equation*}\label{342}
\sum_{j=1}^k \Ks_S(w_1,b_j)\Fs(b_j,w_2)=\delta_{w_1 w_2}.
\end{equation*}
This is true for all $w_1,w_2\in W_S$. Moreover,
$\lim_{|b-w|\rightarrow\infty}\Fs(b,w)=0$, so that by definition
of the inverse real Dirac operator, and by the existence and
uniqueness Lemma \ref{lem13}, we deduce that for all $b\in B_S$,
$w\in W_S$, $\Fs(b,w)=\Ks_S^{-1}(b,w)$.\jump

\noindent Hence, $\mu_n^S(e_1,\ldots,e_k)$ converges to the right
hand side of (\ref{6}) on a subsequence of $n'$s. By Sheffield's
Theorem \cite{Sheffield}, this is the unique limit of the
Boltzmann measures $\mu_n^S$, so that we have convergence for
every $n$.
\end{proof}

\begin{lem}\label{lem11}
{\scriptsize
\begin{equation}\label{7}
\left(\prod_{i=1}^{k} \Ks_S(w_i,b_i)\right)\det_{1 \leq i,\;j \leq
k}\left(\Ks_S^{-1}(b_i,w_j)\right)=\left(\prod_{i=1}^{k}
K_S(w_i,b_i)\right)\det_{1 \leq i,\;j \leq k}
\left(K_S^{-1}(b_i,w_j)\right).
\end{equation}}
\end{lem}
\begin{proof}
By definition of the determinant, the left hand side of (\ref{7})
is equal to
    {\scriptsize
\begin{equation*}\label{416}
\sum_{\sigma\in\Perm_n}\mbox{sgn}\,\sigma\,\left(\prod_{i=1}^{k}
\Ks_S(w_i,b_i)\right)\Ks_S^{-1}(b_1,w_{\sigma(1)})\ldots
\Ks_S^{-1}(b_k,w_{\sigma(k)}),
\end{equation*}}

\noindent where $\Perm_n$ is the set of permutations of $n$
elements. A permutation $\sigma\in\Perm_n$ can be written as a
product of disjoint cycles, so let us treat the case of each cycle
separately. Refer to Section \ref{subsubsec422} for the definition
of the function $\fs_{wx}$.\jump

\noindent $\bullet$ Suppose that in the product there is a
$1$-cycle, that is a point $j$ such that $\sigma(j)=j$. Then,
using Remark \ref{rem1} and Lemma \ref{lem13}, we obtain
\begin{eqnarray*}
\Ks_S(w_j,b_j)&=& \fs_{w_j b_j}(0)K_S(w_j,b_j),\\
\Ks_S^{-1}(b_j,w_j)&=&\overline{\fs_{w_j
b_j}(0)}K_S^{-1}(b_j,w_j).
\end{eqnarray*}
Moreover, $\overline{\fs_{w_j b_j}(0)}=\fs_{w_j b_j}(0)^{-1}$,
hence
\begin{equation}\label{424}
\Ks_S(w_j,b_j)\Ks_S^{-1}(b_j,w_j)=K_S(w_j,b_j)K_S^{-1}(b_j,w_j).
\end{equation}
$\bullet$ Suppose that in the product there is an $\ell$-cycle,
with $\ell\neq 1$. To simplify notations, let us assume
$\sigma(1)=2,\ldots,\sigma(\ell)=1$, and let us prove the
following (indices are written cyclically, i.e. $\ell+1\equiv 1$),
\begin{equation}\label{419}
\prod_{j=1}^\ell
\Ks_S(w_j,b_j)\Ks_S^{-1}(b_j,w_{j+1})=\prod_{j=1}^\ell
K_S(w_j,b_j)K_S^{-1}(b_j,w_{j+1}).
\end{equation}
Again, using Remark \ref{rem1} and Lemma \ref{lem13}, we obtain
    {\small
\begin{equation*}\label{422}
\prod_{j=1}^\ell \Ks_S(w_j,b_j)\Ks_S^{-1}(b_j,w_{j+1})=
\prod_{j=1}^\ell K_S(w_j,b_j)K_S^{-1}(b_j,w_{j+1})\fs_{w_j
b_j}(0)\fs_{w_{j+1} b_j}(0)^{-1}.
\end{equation*}}
Using the definition of the function $\fs_{wx}$, and the fact that
it is well defined yields
    {\scriptsize
\begin{equation*}\label{423}
\prod_{j=1}^\ell\fs_{w_j b_j}(0)\fs_{w_{j+1} b_j}(0)^{-1}=
\prod_{j=1}^\ell\fs_{w_j b_j}(0)\fs_{b_j w_{j+1}}(0)=
\prod_{j=1}^\ell\fs_{w_j w_{j+1}}(0)=\fs_{w_1 w_1}(0)=1.
\end{equation*}}
This proves equation (\ref{419}). Combining equations (\ref{424}),
(\ref{419}), and the fact that every permutation is a product of
cycles, we obtain Lemma \ref{lem11}.
\end{proof}

\begin{lem}\label{lem12}
    {\scriptsize
\begin{equation*}
\left(\prod_{i=1}^{k} K_S(w_i,b_i)\right)\det_{1 \leq i,\;j \leq
k} \left(K_S^{-1}(b_i,w_j)\right)= \left(\prod_{i=1}^{k}
K(w_i,b_i)\right)\det_{1 \leq i,\;j \leq k}
\left(K^{-1}(b_i,w_j)\right).
\end{equation*}}
\end{lem}
\begin{proof}
Since $P$ is simply connected, for every $i,j=1,\ldots,k$, it
contains a path of $\widetilde{R}$ from $w_j$ to $b_i$. Moreover
by Theorem \ref{thm3}, the coefficient of the inverse Dirac
operator corresponding to $b_i,w_j$ only depends on such a path.
Hence $K_{S}^{-1}(b_i,w_j)=K^{-1}(b_i,w_j)$. We also have
$\forall\,i=1,\ldots,k,\, K_S(w_i,b_i)=K(w_i,b_i)$, so that we
deduce Lemma \ref{lem12}.
\end{proof}
\end{proof}

\subsection{Proof of Theorem \ref{thm2}}\label{subsec44}

\noindent The edges of the graph $R^*$ form a countable set. For
every $i \in \NN$, define $f_i:\M(R^*)\rightarrow \{0,1\}$ by
\begin{equation*}\label{45}
f_i(M)=\left\{
\begin{array}{cl}
1 & \mbox{ if the edge $e_i$ belongs to $M$,}\\
0 & \mbox{ else.}
\end{array}
\right.
\end{equation*}
Fix $k\in\NN$, and a $k$-tuple $(s_1,\ldots,s_k)$ of distinct
elements of $\NN$. Let $H\in\B\{0,1\}^k$, where $\B\{0,1\}^k$
denotes the Borel $\sigma$-field of $\{0,1\}^k$, and define a {\bf
cylinder of rank $k$} by
\begin{equation*}\label{46}
A_{(s_1,\ldots,s_k)}(H)=\{ M \in \M(R^*) |
(f_{s_1}(M),\ldots,f_{s_k}(M))\in H\}.
\end{equation*}
Then $A_{(s_1,\ldots,s_k)}(H)$ can be written as a disjoint union
of cylinder sets,
\begin{equation*}\label{47}
A_{(s_1,\ldots,s_k)}(H)=\bigcup_{i=1}^m
\{e_{t_{i1}},\ldots,e_{t_{i \ell_i}}\},
\end{equation*}
(recall that for every $i$, $\{e_{t_{i1}},\ldots,e_{t_{i
\ell_i}}\}$ denotes the set of dimer configurations of $R^*$
containing the edges $e_{t_{i1}},\ldots,e_{t_{i \ell_i}}$). Define
\begin{equation*}\label{48}
\mu_{(s_1,\ldots,s_k)}(H)=\sum_{i=1}^{m}\left(\left(\prod_{j=1}^{\ell_i}
K(w_{t_{ij}},b_{t_{ij}})\right)\det_{1\leq j,k\leq
\ell_i}(K^{-1}(b_{t_{ij}},w_{t_{ik}}))\right).
\end{equation*}

\noindent Let $\Ps$ be a finite simply connected sub-graph of $\Rs$
such that, for every $i=1,\ldots,m$, $P^*$ contains the edges
$e_{t_{i1}},\ldots,e_{t_{i \ell_i}}$. Let $\Ss$ be a periodic
rhombus tiling of the plane that contains $\Ps$ (given by
Proposition \ref{prop1}). Then, by Proposition \ref{prop2},
$\mu_{(s_1,\ldots,s_k)}(H)=\lim_{n\rightarrow\infty}\mu_n^S(A_{(s_1,\ldots,s_k)}(H))$.
From this we deduce that for every $k$, and for every $k$-tuple
$(s_1,\ldots,s_k)$, $\mu_{(s_1,\ldots,s_k)}$ is a probability
measure on $\B\{0,1\}^k$. Moreover, we deduce that the system of
measures $\{\mu_{(s_1,\ldots,s_k)}:$ $(s_1,\ldots,s_k)$ is a
$k$-tuple of distinct elements of $\NN\}$ satisfy Kolmogorov's two
consistency conditions. Applying Kolmogorov's extension theorem,
we obtain the existence of a unique measure $\mu^R$, which
satisfies (\ref{GibbsM}).\jump

\noindent Using the fact that the measure $\mu^R$ of a cylinder
set is the limit of Boltzmann measures, we deduce that the
measure $\mu^R$ is a Gibbs measure in the sense given in the
introduction.\jump

\noindent Assume moreover that the graph $R^*$ is doubly periodic.
Then, for every cylinder set $\{e_1,\ldots,e_k\}$ of $R^*$,
$\mu^R(e_1,\ldots,e_k)=\lim_{n\rightarrow\infty}\mu_n^R(e_1,\ldots,e_k)$.
Moreover by \cite{Sheffield} (see also \cite{KeOS}), the Boltzmann
measures $\mu_n^R$ converge to the minimal free energy per
fundamental domain Gibbs measure which is unique, so that this
proves Theorem \ref{thm2}.
\hspace*{\fill} $\square$ \\

\section{Gibbs measure on the set of all triangular quadri-tilings}\label{sec5}

Recall that $\Q$ is the set of all triangular quadri-tilings up to
isometry, i.e. the set of all quadri-tilings whose underlying
tiling is a lozenge tiling of the equilateral triangular lattice
$\TT$. Denote by $\M$ the set of dimer configurations
corresponding to quadri-tilings of $\Q$. In this section, we first
define the notion of Gibbs measure on $\M$; then we define a
$\sigma$-algebra $\sigma(\B)$ on $\M$, and give an explicit
expression for a Gibbs measure $\mu$ on $(\M,\sigma(\B))$. We
conjecture $\mu$ to be of minimal total free energy per
fundamental domain among a four parameter family of ergodic Gibbs
measures.\jump

\noindent The notion of Gibbs measure on $\Q$ is a natural
extension of the notion of Gibbs measure on dimer configurations
of a fixed graph. Assume a weight function $\nu$ is assigned to
quadri-tiles of triangular quadri-tilings of $\Q$, then a {\bf
Gibbs measure on $\Q$} is a probability measure with the following
property. If a triangular quadri-tiling is fixed in an annular
region, then quadri-tilings inside and outside of the annulus are
independent. Moreover, the probability of any interior triangular
quadri-tiling is proportional to the product of the weights of the
quadri-tiles. Using the bijection between $\Q$ and $\M$, this
yields the definition of a Gibbs measure on $\M$.\jump

\noindent Define $\LL$ to be the set of lozenge-with-diagonals
tilings of the plane, up to isometry. Define $\LL^*$ to be the
graph (which is not planar) obtained by superposing the dual
graphs $L^*$ of lozenge-with-diagonals tilings $L \in \LL$.
Although some edges of $\LL^*$ have length $0$, we think of them
as edges of the one skeleton of the graph, so that to every edge
of $\LL^*$ there corresponds a unique quadri-tile in a
lozenge-with-diagonals tiling of $\LL$. Let $e$ be an edge of
$\LL^*$, and let $q_e$ be the corresponding quadri-tile, then
$q_e$ is made of two adjacent right triangles. If the two
triangles share the hypotenuse edge, they belong to two adjacent
lozenges; else if they share a leg, they belong to the same
lozenge. Let us call these lozenge(s) the {\bf lozenge(s)
associated to the edge $e$}, and denote it/them by $\ls_e$ (that
is $\ls_e$ consists of either one or two lozenges). Let $\ks_e$ be
the edge(s) of $\TT^*$ corresponding to the lozenge(s) $\ls_e$.\\
Let $e_1,\ldots,e_k$ be a subset of edges of $\LL^*$, and define
the {\bf cylinder set} $\{e_1,\ldots,e_k\}$ of $\LL^*$ to be the
set of dimer configurations of $\M$ which contain these edges. Let
us call {\bf connected cylinder} any cylinder of $\LL^*$ which has
the property that the lozenge(s) associated to its edges form a
connected path. Then every cylinder of $\LL^*$ can be expressed as
a disjoint union of connected cylinders. Consider $\B$ the field
consisting of the empty set and of the finite disjoint unions of
connected cylinders. Denote by $\sigma(\B)$ the $\sigma$-field
generated by $\B$.\jump

\noindent Let $\mu^{\TT}$ be the Gibbs measure on dimer
configurations of the honeycomb lattice $\M(\TT^*)$ given in
\cite{Kenyon}. As we have noted before, Theorem \ref{thm2} is true
for all doubly periodic isoradial graphs with critical weights on
their edges. Then, the measure $\mu^{\TT}$ coincides with the
minimal free energy per fundamental domain Gibbs measure given by
Theorem \ref{thm2}, when the doubly periodic isoradial graph is
the honeycomb lattice. As a corollary to Theorem \ref{thm2}, we
obtain:

\begin{corollary}\label{cor4}
There is a probability measure $\mu$ on $(\M,\sigma(B))$ such that
for every connected cylinder $\{e_1,\ldots,e_k\}$ of $\LL^*$,
\begin{equation}\label{430}
\mu(e_1,\ldots,e_k)=\mu^L(e_1,\ldots,e_k)\mu^{\TT}(\ks_{e_1},\ldots,\ks_{e_k}),
\end{equation}
where $L$ is the lozenge-with-diagonals tiling corresponding to
any lozenge tiling $\Ls$ which contains the lozenges
$\ls_{e_1},\ldots,\ls_{e_k}$. Moreover $\mu$ is a Gibbs measure on
$\M$, where the critical weight function is assigned to
quadri-tiles.
\end{corollary}
\begin{proof}
Expression (\ref{430}) is well defined, i.e. independent of the
lozenge tiling $\Ls$ which contains the lozenges
$\ls_{e_1},\ldots,\ls_{e_k}$. Indeed, by definition of a connected
cylinder set, the lozenges associated to the edges
$e_1,\ldots,e_k$ form a connected path of lozenges, say $\gamma$.
Let $\Ls$ be a lozenge tiling that contains $\gamma$, and denote
by $K_{L}$ the complex Dirac operator indexed by the vertices of
the graph $L^*$. Then $K_{L}^{-1}(b_i,w_j)$ is independent of the
lozenge tiling $\Ls$ which contains $\gamma$, indeed
$K_{L}^{-1}(b_i,w_j)$ only depends on an edge-path of
$\widetilde{R}$ from $w_j$ to $b_i$, and since $\Ls$ contains
$\gamma$ which is connected, we can choose the edge-path to be the
same for all such lozenge-with-diagonals tilings $L$. We then use
the fact that $\mu^{\TT}$ and $\mu^{L}$ are probability measures to
prove the two conditions of Kolmogorov's extension theorem.\jump

\noindent The measure $\mu$ is a Gibbs measure on $\M$ as a
consequence of the fact that $\mu^L$ and $\mu^{\TT}$ are Gibbs measure
on $\M(L^*)$ and $\M(\TT^*)$ respectively.
\end{proof}

\noindent Assume that a weight function is assigned to
quadri-tiles of quadri-tilings of $\Q$. Denote by $\Lambda$ the
lattice which acts periodically on $\TT$, and by
$\TT_n=\TT/n\Lambda$, moreover suppose that $\TT_n^*$ is bipartite
(this is possible by eventually replacing $\Lambda$ by
$2\Lambda$). Define $\Q_n$ to be the set of triangular
quadri-tilings whose underlying tiling is a lozenge tiling of
$\TT_n$. Denote by $\mu_n$ the Boltzmann measure on $\Q_n$; that
is the probability of having a given subset of quadri-tiles in a
quadri-tilings of $\Q_n$ chosen with respect to $\mu_n$ is
proportional to the product of the weights of the quadri-tiles. We
make the first following conjecture.

\begin{conjecture}
Suppose the critical weight function is assigned to quadri-tiles,
then the Gibbs measure of Corollary \ref{cor4} is the limit of the
Boltzmann measures $\mu_n$.
\end{conjecture}

\noindent Now fix $(s,t,p,q)\in\RR^4$, and let $\Q_n^{(s,t,p,q)}$
be the subset of quadri-tilings of $\Q_n$ whose first height
change is $(\lfloor ns \rfloor,\lfloor nt\rfloor)$ and second
height change is $(\lfloor np\rfloor,\lfloor nq\rfloor)$. Assuming
that $\Q_n^{(s,t,p,q)}$ is non empty, let $\mu_n^{(s,t,p,q)}$ be
the conditional measure induced by $\mu_n$ on $\Q_n^{(s,t,p,q)}$.
Denote by $Z_n^{(s,t,p,q)}$ the weighted sum of quadri-tilings of
$\Q_n^{(s,t,p,q)}$, and define the {\bf total free energy per
fundamental domain} $\sigma(s,t,p,q)$ by:
\begin{equation*}
\sigma(s,t,p,q)=-\lim_{n\rightarrow\infty}\frac{1}{n^2}\log Z_n^{(s,t,p,q)}.
\end{equation*}
Then the following conjecture is inspired by a result of
\cite{Sheffield}, see also \cite{KeOS}.

\begin{conjecture}
For each $(s,t,p,q)$ for which $\Q_n^{(s,t,p,q)}$ is non empty for
$n$ sufficiently large, $\mu_n^{(s,t,p,q)}$ converges as
$n\rightarrow\infty$ to an ergodic Gibbs measure $\mu^{(s,t,p,q)}$
of slope $(s,t,p,q)$. Furthermore $\mu_n$ itself converges to
$\mu^{(s_0,t_0,p_0,q_0)}$ where $(s_0,t_0,p_0,q_0)$ is the limit
of the slopes of $\mu_n$. If $(s_0,t_0,p_0,q_0)$ lies in the
interior of the set $(s,t,p,q)$ for which $\Q_n^{(s,t,p,q)}$ is
non-empty for $n$ sufficiently large, then every ergodic Gibbs
measure of slope $(s,t,p,q)$ is of the form $\mu^{(s,t,p,q)}$ for
some $(s,t,p,q)$ as above; that is $\mu^{(s,t,p,q)}$ is the unique
ergodic Gibbs measure of slope $(s,t,p,q)$. Moreover, the measure
$\mu^{(s_0,t_0,p_0,q_0)}$ is the unique one which has minimal
total free energy per fundamental domain.
\end{conjecture}
As a consequence, when the critical weight function is assigned to
quadri-tiles, we conjecture the minimal total free energy per
fundamental domain Gibbs measure to be given by the explicit
expression of Corollary \ref{cor4}.

\section{Asymptotics in the case of quadri-tilings}\label{sec6}

Section \ref{subsec61} aims at giving a precise statement of
Theorem \ref{thm1} of the introduction (see Theorem~\ref{thm7}
below), and Section \ref{subsec62} gives consequences of
Theorem~\ref{thm7} for the Gibbs measure $\mu^R$ of Theorem
\ref{thm2} and $\mu$ of Corollary \ref{cor4}.\vspace{-0.1cm}

\subsection{Asymptotics of the inverse Dirac operator}\label{subsec61}

\begin{figure}[h]
\begin{center}
\includegraphics[height=2.5cm]{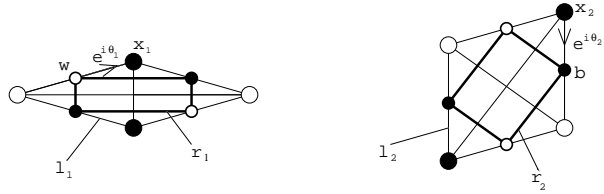}
\end{center}
\caption{Rhombi with diagonals $\ell_1$, $\ell_2$ and their dual
graphs $r_1$, $r_2$.} \label{fig9}
\end{figure}

\noindent Refer to Figure \ref{fig9} for the following notations.
Let $\ell_1'$, $\ell_2'$ be two disjoint side-length two rhombi in
the plane, and let $\ell_1$, $\ell_2$ be the corresponding
rhombi with-diagonals. Assume $\ell_1$ and $\ell_2$ have a fixed
black and white bipartite coloring of their faces. Let $r_1$ and
$r_2$ be the dual graphs of $\ell_1$ and $\ell_2$ ($r_1$ and $r_2$
are rectangles), with the corresponding bipartite coloring of the
vertices. Let $w$ be a white vertex of $r_1$, and $b$ a black
vertex of $r_2$, then $w$ (resp. $b$) belongs to a boundary edge
$e_1$ of $\ell_1$ (resp. $e_2$ of $\ell_2$). By Lemma \ref{lem2},
to the bipartite coloring of the faces of $\ell_1$ and $\ell_2$,
there corresponds a bipartite coloring of the vertices of
$\ell_1'$ and $\ell_2'$. Let $x_1$ (resp. $x_2$) be the black
vertex of the edge $e_1$ (resp. $e_2$). Orient the edge $w x_1$
from $w$ to $x_1$, and let $e^{i \theta_1}$ be the corresponding
vector. Orient the edge $x_2 b$ from $x_2$ to $b$, and let $e^{i
\theta_2}$ be the corresponding vector. Assume $\ell_1'$ and
$\ell_2'$ belong to a rhombus tiling of the plane $\Rs$. Moreover,
suppose that the bipartite coloring of the vertices of ${{R}^*}$
is compatible with the bipartite coloring of the vertices of $r_1$
and $r_2$.
\ \\

\noindent Then we have the following asymptotics for the inverse
Dirac operator $K^{-1}$ indexed by the vertices of ${R}^*$. Refer
to the introduction for comments about Theorem \ref{thm7}.\vspace{0.5cm}

\begin{thm}\label{thm7}
As $|b-w|\rightarrow\infty$, $K^{-1}(b,w)$ is equal to \jump

\noindent{\small $ \frac{1}{2 \pi} \left(
\frac{1}{b-w}+\frac{e^{-i(\theta_1+\theta_2)}}{\bar{b}-\bar{w}}
\right)+\frac{1}{2 \pi}\left( \frac{e^{2 i \theta_1}+e^{2 i
\theta_2}}{(b-w)^3}+\frac{e^{-i(3 \theta_1+\theta_2)} +
e^{-i(\theta_1 + 3\theta_2)}}{(\bar{b}- \bar{w})^3}\right)+O
\left( \frac{1}{|b-w|^3}\right), $}\jump

\noindent where $\theta_1$ and $\theta_2$ are defined above.
\end{thm}

\begin{proof}
\noindent Let us define an edge-path from $w$ to $b$ in
$\widetilde{R}$ (the set of rhombi associated to the edges of
${{R}^*}$). Consider the bipartite coloring of the vertices of
$\Rs$ (given by Lemma \ref{lem2}) associated to the bipartite
coloring of the vertices of ${{R}^*}$. We define the graph $ \Ns$
as follows. Vertices of $ \Ns$ are black vertices of $\Rs$, and
two vertices of $ \Ns$ are connected by an edge if they belong to
the same rhombus in $\Rs$. The graph $\Ns$ is connected because
$\Rs$ is. Each face of $ \Ns$ is inscribable in a circle of radius
two. The circumcenter of a face of $ \Ns$ is the intersection of
the rhombi in $\Rs$, to which the edges on the boundary cycle of
the face belong. Thus the circumcenter is in the closure of the
face, and so faces of $ \Ns$ are convex. Note that the vertices
$x_1$ and $x_2$ are vertices of
the graph $ \Ns$.\jump

\noindent Denote by $(x,y)$ the line segment from a vertex $x$ to a vertex
$y$ of $ \Ns$. An edge $uv$ of $ \Ns$ is called a {\bf
forward-edge} for the segment $(x,y)$ if $<~v-u,y-x>\, \geq 0$. An
edge-path $v_1, \ldots, v_k$ of $ \Ns$ is called a {\bf
forward-path} for the segment $(x,y)$, if all the edges $v_i
v_{i+1}$ are forward-edges for $(x,y)$. Similarly to what has been
done in \cite{Kenyon1}, let us define a forward-path of $ \Ns$ for
the segment $(x_1,x_2)$, from $x_1$ to $x_2$ (see Figure
\ref{fig10}). Let $F_1,\ldots,F_\ell$ be the faces of $ \Ns$ whose
interior intersect $(x_1,x_2)$ (if some edge of $ \Ns$ lies
exactly on $(x_1,x_2)$, perturb the segment $(x_1,x_2)$ slightly,
using instead a segment $(x_1 + \varepsilon_1,x_2+\varepsilon_2)$
for two generic infinitesimal translations
$(\varepsilon_1,\varepsilon_2$)). Note that the number of such
faces is finite because the rhombus tiling of the plane $\Rs$ has
only finitely many different rhombi. Then for $j=1,\ldots,\ell-1$,
$F_j \cap F_{j+1}$ is an edge $e_{j+1}$ of $ \Ns$ crossing
$(x_1,x_2)$. Set $v_1=x_1$, $v_\ell=x_2$, and for
$j=1,\ldots,\ell-2,$ define $v_{j+1}$ to be the vertex of
$e_{j+1}$ such that the edge $e_{j+1}$ oriented towards $v_{j+1}$
is a forward-edge for $(x_1,x_2)$. Then, for $j=1,\ldots,\ell-1$,
the vertices $v_j$ and $v_{j+1}$ belong to the face $F_j$. Take an
edge-path from $v_j$ to $v_{j+1}$ on the boundary cycle of $F_j$,
such that it is a forward-path for $(x_1,x_2)$. Such a path exists
because faces of $ \Ns$ are convex. Thus, we have built a
forward-path of $ \Ns$ for $(x_1,x_2)$, from $x_1$ to $x_2$.
Denote by $u_1=x_1,u_2,\ldots,u_{k-1},u_k=x_2$ the vertices of
this path.

\begin{figure}[h]
\begin{center}
\includegraphics[height=5.6cm]{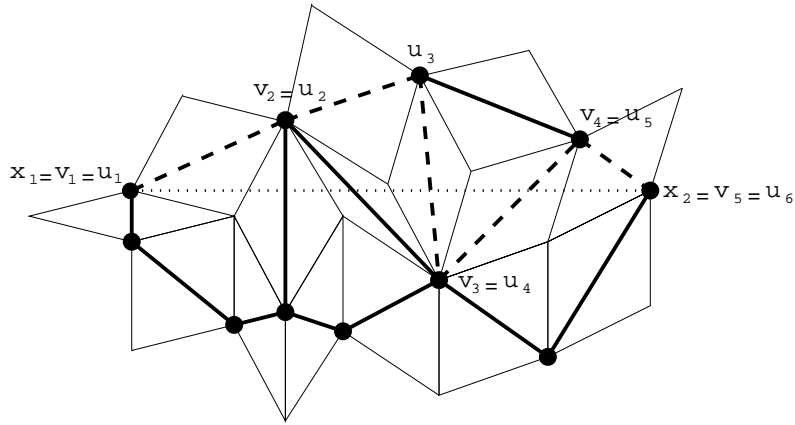}
\end{center}
\caption{forward-path from $x_1$ to $x_2$ for the segment
$(x_1,x_2)$.} \label{fig10}
\end{figure}

\noindent Let us now define an edge-path of $\widetilde{R}$ from
$w$ to $b$. Note that the edges $w x_1$ and $x_2 b$ are edges of
$\widetilde{R}$. For $j=1,\ldots,k-1,$ define the following
edge-path of $\widetilde{R}$ from $u_j$ to $u_{j+1}$ (see Figure
\ref{fig11}). Remember that $u_j u_{j+1}$ is the diagonal of a
rhombus of $\Rs$, say $\tilde{\ell_j'}$. Let $\tilde{r_j}$ be the
dual graph of $\tilde{\ell_j}$. Let $u_j^1$ be the black vertex in
$\tilde{r_j}$ adjacent to $u_j$, let $u_j^2$ be the crossing of
the diagonals of $\tilde{\ell_j}$, and let $u_j^3$ be the white
vertex in $\tilde{r_j}$ adjacent to $u_{j+1}$. Then the path $u_j,
u_j^1,u_j^2,u_j^3,u_{j+1}$ is an edge-path of $\widetilde{R}$.
Thus
$w,x_1=u_1,u_1^1,u_1^2,u_1^3,u_2,\ldots,u_{k-1},u_{k-1}^1,u_{k-1}^2,u_{k-1}^3,u_k=x_2,b$
is an edge-path of $\widetilde{R}$, from $w$ to $b$. Orient the
edges in the path towards the black vertices of ${{R}^*}$, and
away from the white vertices of ${{R}^*}$.

\begin{figure}[h]
\begin{center}
\includegraphics[height=2.1cm]{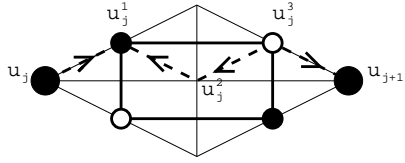}
\end{center}
\caption{Edge-path of $\widetilde{R}$ from $u_j$ to $u_{j+1}$.}
\label{fig11}
\end{figure}

\noindent Let $e^{i \beta_j^1}$, $e^{i \beta_j^2}$, $e^{i
\alpha_j^1}$, $e^{i \alpha_j^2}$ be the vectors corresponding
respectively to the edges $u_j u_j^1$, $u_j^3 u_{j+1}$, $u_j^2
u_j^1$, $u_j^3 u_j^2$. Without loss of generality, suppose that
$x_2-x_1$ is real and positive. Then for $j=1,\ldots,k-1$, and
$\ell=1,2$, we have:
$$
\cos \beta_j^\ell - \cos \alpha_j^\ell = \frac{<u_{j+1}-u_j, x_2-x_1>}{2|x_2-x_1|}.
$$
Since $u_1,\ldots,u_k$ is a forward-path for $(x_1,x_2)$, this
quantity is positive, thus $\cos \beta_j^\ell \geq \cos \alpha_j^\ell$.

\noindent Moreover, since there is only a finite number of
different rhombi in $\Rs$, $k=O(|b-w|)$. For the same reason,
there is a finite number of angles $\beta_j^\ell$, and they are
all in $[-\pi+\Delta, \pi - \Delta]$, for some small $\Delta
> 0$ (in the general case where the angle of the vector $x_2 -x_1$
is $\theta_0$, the angles $\beta_j^\ell$ would be in the interval
$[\theta_0-\pi+\Delta, \theta_0 + \pi - \Delta]$). Thus by Theorem
$4.3$ of \cite{Kenyon1}, we have that $K^{-1}(b,w)$ is equal to:
    {\small
\begin{equation}\label{thmkenyon}
\frac{1}{2 \pi} \left(
\frac{1}{b-w}+\frac{\gamma}{\bar{b}-\bar{w}}\right)
    +\frac{1}{2 \pi} \left( \frac{\xi_2}{(b-w)^3}+
    \frac{\gamma \bar{\xi_2}}{(\bar{b}- \bar{w})^3}\right)
    +O \left( \frac{1}{|b-w|^3}\right),
\end{equation}}
where $\gamma=e^{-i(\theta_1+\theta_2)} \displaystyle
\prod_{j=1}^{k-1} \prod_{\ell=1}^{2} e^{i(-\beta_j^\ell+\alpha_j^\ell)}$,
and $\xi_2=e^{2i\theta_1}+e^{2i\theta_2} + \displaystyle
\sum_{j=1}^{k-1} \sum_{\ell=1}^{2} e^{2i\beta_j^\ell} -
e^{2i\alpha_j^\ell}$.\\

\noindent Note that for $j=1,\ldots,k-1$, we have $\alpha_j^2
\equiv (\beta_j^1+\pi) \;{\rm mod} [2\pi]$, and $\beta_j^2 \equiv
(\alpha_j^1 + \pi) \;{\rm mod} [2\pi]$, thus:
\begin{eqnarray*}
\prod_{\ell=1}^{2} e^{i(-\beta_j^\ell+\alpha_j^\ell)} &=&
    e^{i(-\beta_j^1+\alpha_j^1)}e^{i(-\alpha_j^1-\pi+\beta_j^1+\pi)}=1,\\
\sum_{\ell=1}^{2} e^{2i\beta_j^\ell} - e^{2i\alpha_j^\ell} &=&
    e^{2i\beta_j^1} - e^{2i\alpha_j^1}+e^{2i(\alpha_j^1+\pi)}-e^{2i(\beta_j^1 + \pi)}=0.
\end{eqnarray*}
Therefore $\gamma=e^{-i(\theta_1+\theta_2)}$,
$\xi_2=e^{2i\theta_1}+e^{2i\theta_2}$, which proves the theorem.
\end{proof}

\subsection{Asymptotics of the Gibbs measures on quadri-tilings}\label{subsec62}

\noindent Let $\Rs$ be a rhombus tiling of the plane, and $R$ be
the corresponding rhombus-with-diagonals tiling. Consider a subset
of edges $e_1=w_1 b_1,\ldots,e_k=w_k b_k$ of $R^*$, and recall
that $\mu^R$ is the Gibbs measure on $\M(R^*)$ given by Theorem
\ref{thm2}.

\begin{corollary}\label{cor1}
When $\forall\,j \neq i$, $|w_{j}-b_{i}| \rightarrow \infty$, then
up to the second order term, $\mu^R(e_1,\ldots,e_k)$ only depends
on the rhombi of $\Rs$ to which the vertices $b_i$ and $w_j$
belong, and else is independent of the structure of the graph
$\Rs$.
\end{corollary}

\begin{proof}
This is a consequence of the explicit formula for $\mu^R(e_1,\ldots,e_k)$ of
Theorem \ref{thm2}, and of the asymptotic formula for the inverse
Dirac operator of Theorem \ref{thm7}.
\end{proof}

\noindent Recall that $\LL^*$ is the non-planar graph obtained by
superposing duals of lozenge-with-diagonals tilings of $\LL$. Let
$e_1=w_1 b_1,\ldots,e_k=w_k b_k$ be a subset of edges of $\LL^*$.
Define $\LL^E$ to be the set of lozenge-with-diagonals tilings of
the plane that contain the lozenges associated to the edges
$e_1,\ldots,e_k$.

\begin{corollary}\label{cor2}
When $\forall\,j \neq i$, $|w_{j}-b_{i}| \rightarrow \infty$, then
up the second order term, $\mu^L(e_1,\ldots,e_k)$ is independent
of $L \in \LL^E$.
\end{corollary}

\begin{proof}
As in Section \ref{sec5}, we choose an embedding of $\LL^*$ so
that every edge of $\LL^*$ uniquely determines the lozenge(s) it
belongs to. Corollary \ref{cor2} is then a restatement of
Corollary \ref{cor1}.
\end{proof}

\noindent Recall that $\mu$ is the Gibbs on triangular
quadri-tilings given by Corollary \ref{cor4}. Let
$\ls_{e_1},\ldots,\ls_{e_k}$ be the lozenges associated to the
edges $e_1,\ldots,e_k$ of $\LL^*$, and let
$\ks_{e_1},\ldots,\ks_{e_k}$ be the corresponding edges of
$\TT^*$.

\begin{corollary}\label{cor3}
When $\forall\,j \neq i$, $|w_{j}-b_{i}|\rightarrow \infty$, and
for every $L\in \LL^E$, we have
\begin{equation*}
\mu(e_1,\ldots,e_k)=
\mu^L(e_1,\ldots,e_k)\mu^{\TT}(\ls_{e_1},\ldots,\ls_{e_k})+O\left(\frac{1}{(\bar{b}-\bar{w})^3}\right).
\end{equation*}
\end{corollary}

\begin{proof}
This is a consequence of the explicit formula for
$\mu(e_1,\ldots,e_k)$ of Corollary \ref{cor4}, and of Corollary
\ref{cor2}.
\end{proof}

\end{document}